\documentclass[10pt]{article}

\usepackage{geometry}

\usepackage{amsmath}
\usepackage{amssymb}
\usepackage{graphicx}
\usepackage{amsthm}
\usepackage{caption}
\usepackage{tikz-cd} 
\usepackage{xcolor}
\usepackage{scrextend}
\usepackage{mathtools}
\usepackage{xcolor}
\usepackage{bbm}
\usepackage{tabu}
\usepackage{tcolorbox}
\usepackage{dsfont}
\usepackage{pgfplots}
\usetikzlibrary{svg.path}
\usepackage{amsthm}
\usepackage{accents}
\usepackage{fancyhdr}
\usepackage{changepage}
\usepackage[hang,flushmargin]{footmisc}
 \usepackage[normalem]{ulem}
 \usepackage{color,soul}
 \usepackage{comment}
 
\usepackage{listings}
\definecolor{codegreen}{rgb}{0,0.6,0}
\definecolor{codegray}{rgb}{0.5,0.5,0.5}
\definecolor{codepurple}{rgb}{0.58,0,0.82}
\definecolor{backcolour}{rgb}{0.95,0.95,0.92}
\lstdefinestyle{mystyle}{
    backgroundcolor=\color{backcolour},   
    commentstyle=\color{codegreen},
    keywordstyle=\color{magenta},
    numberstyle=\tiny\color{codegray},
    stringstyle=\color{codepurple},
    basicstyle=\ttfamily\footnotesize,
    breakatwhitespace=false,         
    breaklines=true,                 
    captionpos=b,                    
    keepspaces=true,                 
    numbers=left,                    
    numbersep=5pt,                  
    showspaces=false,                
    showstringspaces=false,
    showtabs=false,                  
    tabsize=2
}
\usepackage{cite}

\usepackage{float}
\usepackage{youngtab}

\usepackage[font=small,labelfont=bf]{caption}
\DeclareCaptionFont{small}{\small}

\usepackage[pdftex,pdfpagelabels,bookmarks,hyperindex,hyperfigures]{hyperref}

\makeatletter
\DeclareRobustCommand\widecheck[1]{{\mathpalette\@widecheck{#1}}}
\def\@widecheck#1#2{%
    \setbox\z@\hbox{\m@th$#1#2$}%
    \setbox\tw@\hbox{\m@th$#1%
       \widehat{%
          \vrule\@width\z@\@height\ht\z@
          \vrule\@height\z@\@width\wd\z@}$}%
    \dp\tw@-\ht\z@
    \@tempdima\ht\z@ \advance\@tempdima2\ht\tw@ \divide\@tempdima\thr@@
    \setbox\tw@\hbox{%
       \raise\@tempdima\hbox{\scalebox{1}[-1]{\lower\@tempdima\box
\tw@}}}%
    {\ooalign{\box\tw@ \cr \box\z@}}}
\makeatother

\newcommand{\Z}{\mathbb{Z}}
\newcommand{\N}{\mathbb{N}}
\newcommand{\intersection}{\cap}

\newcommand{\sdp}{\rtimes}

\newcommand{\isom}{\cong}

\newcommand{\inj}{\hookrightarrow}
\newcommand{\quotient}[2]{{\raisebox{.2em}{$#1$}\left/\raisebox{-.2em}{$#2$}\right.}}

\DeclarePairedDelimiterX{\inner}[2]{\langle}{\rangle}{#1, #2}

\newcommand{\C}{\mathbb{C}}

\newcommand{\Ds}{\bigoplus}

\renewcommand{\bar}{\overline}

\newcommand{\Q}{\mathbb{Q}}

\newcommand{\x}{\times}

\renewcommand{\x}{\otimes}

\newcommand{\M}{\mathbb{M}}

\renewcommand{\choose}[2]{\binom{#1}{#2}}

\makeatletter
\newtheorem*{rep@theorem}{\rep@title}
\newcommand{\newreptheorem}[2]{%
\newenvironment{rep#1}[1]{%
 \def\rep@title{#2 \ref{##1}}%
 \begin{rep@theorem}}%
 {\end{rep@theorem}}}
\makeatother

\newtheorem{theorem}{Theorem}[section]
\newreptheorem{theorem}{Theorem}
\newtheorem{lemma}[theorem]{Lemma}
\newreptheorem{lemma}{Lemma}
\newtheorem{property}[theorem]{Property}

\newtheorem{proposition}[theorem]{Proposition}
\newreptheorem{proposition}{Proposition}
\newtheorem{corollary}[theorem]{Corollary}
\newreptheorem{corollary}{Corollary}

\newtheorem*{theorem*}{Theorem}

\newtheorem{definition}[theorem]{Definition}

\newtheorem{observation}[theorem]{Observation}

\renewenvironment{proof}{{\flushleft \emph{\underline{Proof:}}}}{\hfill $\square$}

\DeclareMathOperator{\sep}{ | }

\DeclareMathOperator{\Id}{Id}

\DeclareMathOperator{\Hom}{Hom}

\DeclareMathOperator{\Res}{Res}
\DeclareMathOperator{\Ind}{Ind}

\DeclareMathOperator{\FI}{FI}

\definecolor{myblue}{rgb}{0.70,0.81,1}
\definecolor{myorange}{rgb}{1,0.72,0.4}
\definecolor{myred}{rgb}{1,0.58,0.58}
\definecolor{myyellow}{rgb}{1,0.88,0.56}
\definecolor{mylime}{rgb}{0.9,1,0}
\definecolor{mycyan}{rgb}{0.78,0.94,1}
\definecolor{mypink}{rgb}{0.78,0.94,1}
\definecolor{mymagenta}{rgb}{0.9,0.70,0.87}
\definecolor{myolive}{rgb}{0.83,0.87,0.54}
\definecolor{mypink}{rgb}{0.99,0.74,0.94}
\definecolor{mybrown}{rgb}{0.89,0.68,0.68}
 
 \newtcolorbox{custombox}[3][]
{
  colframe = #2!25,
  colback  = #2!10,
  coltitle = #2!20!black,  
  title    = #3,
  #1,
}

\usepackage{tikz}

\renewcommand{\restriction}{\mathord{\upharpoonright}}

\providecommand{\keywords}[1]
{
  \small	
  \emph{Keywords:} #1
}

\providecommand{\amsclassification}[1]
{
  \small	
  \emph{MSC:} #1
}

\providecommand{\cis}[1]
{
  \small	
  \emph{Competing Interest Statement:} #1
}

\begin{document}

\title{Polynomial relations between operators on chains of representation rings}
\author{Sun Woo Park$^a$, Maithreya Sitaraman$^{b}$}
\date{\small \flushleft $^a$ Department of Mathematics, University of Wisconsin-Madison, 480 Lincoln Dr, Madison, WI, United States \\
    $^b$ Department of Mathematics, Columbia University, 2990 Broadway, New York, NY, United States \\ }
\maketitle
\thispagestyle{empty}

\begin{abstract}
Given a chain of groups $G_0 \le G_1 \le G_2 ... $, we may form the corresponding chain of their representation rings, together with induction and restriction operators. Let $\textrm{Res}^l$ denote the operator which restricts down $l$ steps, and similarly for $\textrm{Ind}^l$. Observe that $\textrm{Ind}^l \textrm{Res}^l$ is an operator from any particular representation ring to itself. We provide explicit rigid constraints that a group of chain with surjective restriction operators satisfying the polynomial property must obey.

The central question that this paper addresses is: ``What happens if the $\textrm{Ind}^l \textrm{Res}^l$ operator is a polynomial in the $\textrm{Ind } \textrm{Res}$ operator?''. It is well known that chains of wreath products $\{H^n \rtimes S_n\}_{n \in \N}$ have this property. In this paper, we deduce rigid numerical constraints any arbitrary chains of groups with surjective restriction operator and the aforementioned polynomial property must satisfy. These polynomials can be obtained from the orders of the groups, and are uniquely parametrized by two integers. We deduce the desired parametrization from deducing character theoretic properties of group chains with the desired polynomials.
\end{abstract}

\amsclassification{20C05 $\cdot$ 20C30}


\cis{Declarations of interest: none.}

\let\thefootnote\relax\footnotetext{Email addresses: spark483@math.wisc.edu (Sun Woo Park), maithreya@math.columbia.edu (Maithreya Sitaraman).}

\flushleft

{\small
\tableofcontents
}

\section{Introduction}

\flushleft

\ 

There has been an ever-growing importance placed on induction and restriction operators, both in the context of symmetric groups as well as in the more general setting of chains of algebraic objects. Recent research focuses on the branching graph perspective for symmetric group representation theory in which induction and restriction operators play the central role \cite{bg3} \cite{bg1} \cite{bg2}. The idea of Gelfand Tsetlin bases has been extended to a variety of other algebraic structures as well, see for example \cite{gz1} \cite{gz2}. Structures that emerge when studying induction and restriction include PSH algebra structures \cite{ir1} and Heisenberg algebra structures \cite{ir2}. Wreath products have naturally arisen when studying Heisenberg algebra structures \cite{wp1} \cite{wp2}, and they respect Littlewood-Richardson type rules \cite{lr1}. 

\

Let $G$ be any finite group, and let $R(G)$ be the vector space of representations of $G$ over an algebraically closed field $k$, whose basis is given by the set of irreducible representations of $G$. Given a chain of groups $G_0 \leq G_1 \leq G_2 \leq \cdots$, we form the corresponding chain of vector spaces of their representations, together with induction and restriction operators.
\begin{center}
\begin{tikzcd}
R(G_0) \arrow[shift left=1.5 ex]{r}{\Ind} & \arrow[shift left=1.5 ex]{l}{\Res}  R(G_1) \arrow[shift left=1.5 ex]{r}{\Ind} & \arrow[shift left=1.5 ex]{l}{\Res} R(G_2) \arrow[shift left=1.5 ex]{r}{\Ind} & \arrow[shift left=1.5 ex]{l}{\Res} ... 
\end{tikzcd}
\end{center}
\ 

Note that $\Ind$ and $\Res$ are considered as maps of vector spaces rather than maps of rings. This is because restrictions respect tensor products, whereas inductions do not. Denote by $\Res^l$ the operator which restricts representations from $G_n$ to $G_{n-l}$, and denote by $\Ind^l$ the respective operator. 

\

In this paper, we pursue the following natural structural question: ``What happens if the $\Ind^l \Res^l$ operator is a polynomial in the $\Ind \Res$ operator over the representation ring $R(G_n)$ of $G_n$?''. In the special case where the chain of group rings form a tower of algebras (i.e. there is an injection $\mathbb{C}[G_k] \otimes \mathbb{C}[G_l] \to \mathbb{C}[G_{k+l}]$), then it is known that such a polynomial property will hold, see for example \cite{dg2} \cite{dg3} \cite{dg1}. In \cite{dg2}, the authors consider towers of algebras satisfying certain axioms, and prove that their Grothendieck rings have dual graded graph structures. They then appeal to \cite{dg1}, which, given a dual graded graph structure, uses combinatorial arguments to express the $\Ind^l \Res^l$ operator as a polynomial in the $\Ind \Res$ operator. In this paper, we pursue this question from a representation theoretic perspective, examine the representation theoretic consequences of the above polynomial property, and explore necessary conditions for such a polynomial property to hold. This yields results of a very different flavor from \cite{dg2} \cite{dg3} and \cite{dg1}. 

\

Throughout this paper, we denote by Property $(*)$ the constraint that the operator $\Ind^l \Res^l$ is a polynomial in terms of the $\Ind \Res$ operator for any $l$.

\

\begin{property}[] \label{*1}
$(*)$ For every $l \in \N$ there exists a polynomial $f_l$ such that $\Ind^l \Res^l = f_l(\Ind \Res)$.
\end{property}

\

This paper is devoted to finding rigid numerical constraints on the chains of groups $\{G_i\}$ that satisfy Property $(*)$. The main result of this paper shows that representation theoretic perspective allows us to rigidly determine the polynomials $f_l$ by using two integral parameters $B$ and $C$. 

\

\begin{theorem}[Polynomials are determined by two parameters] \label{pdtp}
Let $\{G_n\}_{n \ge 0}$ be a non-constant chain of groups with surjective restriction operators ( \hyperref[sc]{surjective chain} of groups) which satisfies \hyperref[*1]{Property $(*)$}; that is, there exists polynomials $\{f_l\}_{l \ge 1}$ such that $\Ind^l \Res^l = f_l(\Ind \Res)$.  Let $a_n = \frac{|G_n|}{|G_{n-1}|}$. Then: \\
\begin{addmargin}[2em]{0em}
(1) $\{a_n\}_{n \in \N}$ is an infinite set. Moreover, there exist two parameters $B \in \N$ and $C \in \Z$ such that for all $n$, $a_n = B a_{n-1} + C$. \\
(2) The two parameters $B$ and $C$ from Part (2) completely determine the polynomials $\{f_l\}_{l \in \N}$. The polynomial $f_l$ can be expressed in terms of these parameters as: $$f_l(X) = \frac{1}{B^{\frac{l(l-1)}{2}}} X (X - C) (X - [1+ B]C) (X - [1 + B + B^2]C) ... (X - [1 + B +  .... +B^{l-2} ] C) $$
\end{addmargin}
\end{theorem}

\

For instance, all polynomials are uniquely determined by the polynomial of $\Ind^2 \Res^2$ (i.e $f_2$), and that $f_2$ takes the form $f_2(X) = \frac{1}{B} X(X-C)$ where $B \in \N$ and $C \in \Z$. Furthermore, $B$ and $C$ can be easily read off from the orders of the groups $\{G_i\}$. 

\

While previous research such as \cite{ir2} and \cite{ir1} present algebraic structures that induction and restriction operators respect, our theorem, on the other hand, explicitly demonstrates rigid constraints which prevent induction and restriction operators from respecting certain algebraic structures. It can be viewed as a partial answer to a generalization of a Conjecture by Christian Gaetz as follows: The collection of chains of groups satisfying property $(*)$ includes that of dual towers of groups, a family of group chains which satisfy the \hyperref[hap]{Heisenberg algebra property}. Gaetz \cite{dg4} conjectured that the numerical constraints the desired polynomial relation satisfy forces the chain of groups to be a chain of wreath products. He proves the conjecture in the special case when the scaling factor is either $1$ or is a prime. Since the \hyperref[hap]{Heisenberg algebra property} implies \hyperref[*1]{Property $(*)$} but not vice versa, our main result can be viewed as a partial answer to a generalization of the conjecture by Gaetz.

\ 

It should be noted that a corollary of our work is that, in the language of dual graded graphs (see  \cite{dg1}), $\Ind$ and $\Res$ operators have a recurrent commutation relation if and only if they have linear commutation relation.  In \cite{dg2} and \cite{dg3}, the authors assume conditions that imply a linear commutation relation and deduce a recurrent commutation relation. However, in the context of ``surjective chains of groups'', we are able to prove the converse (the harder direction). This is a novel aspect of our paper which is not addressed in the literature.

\

\section{Setup and notation}

\

We work with chains of groups $G_0 \le G_1 \le G_2 ... $ and their corresponding chains of representation rings over $\C$:

\begin{center}
\begin{tikzcd}
R(G_0) \arrow[shift left=1.5 ex]{r}{\Ind} & \arrow[shift left=1.5 ex]{l}{\Res}  R(G_1) \arrow[shift left=1.5 ex]{r}{\Ind} & \arrow[shift left=1.5 ex]{l}{\Res} R(G_2) \arrow[shift left=1.5 ex]{r}{\Ind} & \arrow[shift left=1.5 ex]{l}{\Res} ... 
\end{tikzcd}
\end{center}

\

that are surjective, in the sense that 

\begin{definition}[Surjective Chain] \label{sc} A \underline{surjective chain} of groups $\{G_n\}_{n \ge 0}$ is a sequence of groups such that: \\ 
\begin{addmargin}[2em]{0em}
(1) $G_{n-1} \le G_n$ for all $n$ \\
(2) $\Res: R(G_n) \to R(G_{n-1})$ is a surjective linear map for all $n$ \\ 
\end{addmargin}
\end{definition}

\

\begin{adjustwidth}{1cm}{}
\textbf{Remark:} The condition that $\Res: R(G_n) \to R(G_{n-1})$ is a surjection is equivalent by Frobenius reciprocity to the condition that $\Ind: R(G_{n-1}) \to R(G_n)$ is an injection. To see this, simply note that: $\Res$ is not surjective $\iff$ there exists some nonzero representation $w \in R(G_{n-1})$ which is perpendicular to the image of $\Res$ $\iff$ there exists a basis $\{x_1,..,x_n\}$ of $R(G_n)$ such that $\inner{\Res(x_i)}{w} = 0$ for all $i$ $\iff$ there exists a basis $\{x_1,..,x_n\}$ of $R(G_n)$ such that $\inner{x_i}{\Ind(w)} = 0$ for all $i$ $\iff \Ind(w) = 0$ for some nonzero $w$ $\iff$ $\Ind$ is not injective.
\end{adjustwidth}

\ 

\begin{adjustwidth}{1cm}{}
\textbf{Remark:}\label{dual}
We note that the construction of chains of representation rings of $\{G_n\}_{n \geq 0}$ is analogous to the construction of the Grothendieck ring of the tower of algebras $\bigoplus_{n \geq 0} \mathbb{C} G_n$. The only difference is that the Grothendieck ring is a $\Z$-module, while we consider the representation ring as a $\mathbb{C}$-vector space. We refer to Section 2 of \cite{dg2} and Section 3 of \cite{dg3} for the conditions on the tower of algebras $\bigoplus_{n \geq 0} \mathbb{C} G_n$. However, we also note that the conditions on surjective chains of groups $\{G_n\}_{n \geq 0}$ are different from those on the tower of algebras provided in \cite{dg2} and \cite{dg3}. One of the axioms \cite{dg2} and \cite{dg3} assume is that the external multiplication $\mathbb{C} G_m \otimes \mathbb{C} G_n \to \mathbb{C} G_{m+n}$ is an injection. This is a rigid condition on the level of algebras. However, the only condition we have on the level of algebras is the fact that $\mathbb{C}[G_n]$ is a subalgebra of $\mathbb{C}[G_{n+1}]$, which is a weak condition. We do, however, impose a condition on the level of representation rings (rather than the level of algebras) that the restriction operator is surjective. In \cite{dg2} and  \cite{dg3}, Mackey's theorem implies the existence of a Hopf algebra structure on the level of representation rings with multiplication $m: \overline{R} \otimes \overline{R} \to \overline{R}$ given by $\Ind_{G_m \times G_n}^{G_{m+n}}$ and comultiplication $\Delta: \overline{R} \to \overline{R} \otimes \overline{R}$ given by $\sum_{k+l = n} \Res_{G_k \times G_l}^{G_{k+l}}$. Our surjectivity condition need not imply such a Hopf algebra structure.
\end{adjustwidth}

\

Let $\bar{R}\{l\}$ denote the subchain of the chain of representation rings comprising of all components of degree $\ge l$, i.e $\bar{R}\{l\} = \Ds_{i=l}^{\infty} R(G_i)$. We now consider the family of grading-preserving linear operators $\{\Ind^l \Res^l: \bar{R}\{l\} \to \bar{R}\{l\} \}_{l \in \N}$, where $\Ind^l \Res^l: R(G_n) \to R(G_n)$ is the operator which restricts a representation $l$ steps down and then induces the resulting expression $l$ steps up. That is: 

\ 

$$\Ind^l \Res^l \restriction_{R(G_n)} := \begin{cases} \Ind^{G_n}_{G_{n-l}} \Res^{G_n}_{G_{n-l}} & \textrm{ if $n \ge l$ } \\ \textrm{ not defined } & \textrm{ if $n < l$ } \end{cases}$$ 

\ 

This construction allows us to understand the overarching question of this paper: ``What happens if the $\Ind^l \Res^l$ operator is a polynomial in the $\Ind \Res$ operator?''. There exists a polynomial over $\Q$ (in fact, over $\Z$ as seen from Section 4) depending on $l$, $f_l(X) = c_0 + c_1 X + c_2 X^2 + ... + c_{d} X^d$ such that for any $n \ge l$ and any representation $u \in R(G_n)$, 

$$\Ind^l \Res^l(u) = c_0 u + c_1 \Ind \Res(u) + c_2 \Ind \Res \Ind \Res (u) + ... + c_d (\Ind \Res)^d (u)$$

where multiplication on the right hand side is of course composition. In addition, $\Ind^l \Res^l \restriction_{R(G_n)}$ and $\Ind \Res \restriction_{R(G_n)}$ are matrices for every $n \ge l$, and the property above states that there exists a single polynomial $f_l$ such that for all $n \ge l$, the matrix $\Ind^l \Res^l \restriction_{R(G_n)}$ is equal to $f_l$ applied to the matrix $\Ind \Res \restriction_{R(G_n)}$. 

\

The operators $\Ind$ and $\Res$ are considered as maps of vector spaces rather than maps of rings, the reason being that although restrictions respect tensor products, inductions do not. In this paper, we are interested in obtaining information from restrictions, and so it is of the utmost importance that the restriction maps above are surjective. The restriction from the symmetric group to the alternating group $\Res: R(S_n) \to R(A_n)$ is an example of a restriction which is not surjective, and we wish to exclude such cases. Since we think of operators like $\Ind \Res$ as linear maps, we like to denote them by capital letters. Since we think of representations as vectors, we denote them by small letters. The following table of notation outlines the notation used in this paper:

\

\begin{table}[!htb] 
\centering
\begin{tabular}{ |p{5cm}||p{9.5 cm}| }
 \hline
 \textbf{Notation} &  \textbf{What it denotes} \\
  \hline
Lower case letters such as $u$, $w$, $x$, $y$ & Representations living in some $R(G_n)$ \\
    \hline
Lower case letters with a tilde such as $\tilde{w}$ & If $w \in R(G_k)$, $\tilde{w}$ is a lift of $w$ by Restriction that lives in $R(G_n)$ for some $n > k$. That is, $\Res^{n-k}(\tilde{w}) = w$ \\
    \hline
 $\Ind^k$ & The $k$-step induction operator  $\Ind_{G_{n}}^{G_{n+k}}$ for all $n$ \newline \\ 
  \hline
  $\Res^k$ & The $k$-step restriction operator $\Res_{G_{n}}^{G_{n+k}}$ for all $n$ \newline \\ 
  \hline
$X$ & The operator $\Ind \Res$ \\
  \hline
$\delta_h$ & The column of the character table of $G_n$ associated with the conjugacy class $[h]$. As a vector in $R(G_n)$, $\delta_h = \sum_{\textrm{irreps $u$}} \chi_u(h) u$ \\
 \hline
$t_n$ or $t$ & The trivial representation of the group $G_n$ \\
 \hline
 $s_n$ or $s$ & The sign representation of the symmetric group $S_n$ \\
 \hline
 $v_n$ or $v$ & The standard representation of the symmetric group $S_n$ \\
 \hline
 $\wedge^2_n$ or $\wedge^2$ & The representation of the symmetric group $S_n$ corresponding to the young diagram $(n-2,1,1)$ \\
 \hline
 $p_n$ or $p$ & The representation of the symmetric group $S_n$ corresponding to the young diagram $(n-2,2)$ \\
 \hline
 $b_n$ or $b$ & The representation of the symmetric group $S_n$ corresponding to the young diagram $(n-3,3)$. \\
 \hline
 $r_n$ or $r$ & The representation of the symmetric group $S_n$ corresponding to the young diagram $(n-3,2,1)$.\\
 \hline
 $\M(G,u)$ & The McKay graph associated to a group $G$ and its representation $u$ \\
  \hline
  $[h]_n$, $[\tau]_n$ & Conjugacy classes of $h$ and $\tau$ when considered as an element of $G_n$ \\
  \hline
\end{tabular}
\caption{Table of notation}
\label{notations}
\end{table}

\

\section{Motivating examples}

\

In this section we will introduce the motivating examples of chains which satisfy Property $(*)$, namely symmetric groups and (more generally) wreath products. Both of these chains have the additional Heisenberg-algebra property as defined below: 

\

\begin{property}[Heisenberg algebra property] \label{hap} A chain of groups $\{G_n\}_{n \ge 0}$ is said to satisfy the Heisenberg algebra property with scaling $M \in \N$ if: 
$$\Res \Ind  - \Ind \Res  = M \Id$$
\end{property}

\

It is well known that:

\begin{proposition} \label{cwpsc}
Let $H$ be any finite group. 
\begin{enumerate}
    \item The chain of wreath products $\{H^n \sdp S_n\}_{n \in \N}$ is a surjective chain. 
    \item The chain of wreath products $\{H^n \sdp S_n\}_{n \in \N}$ satisfies the Heisenberg algebra property with scaling $|H|$.
\end{enumerate}
\end{proposition}

\

In particular, setting $H = \{1\}$, the chain of symmetric groups satisfies the Heisenberg algebra property with $M=1$. It remains to show that Heisenberg algebra property implies \hyperref[*1]{Property $(*)$}. This is also a well known fact; for instance, the reader may see for example Corollary 1.4.11 of \cite{dg1}. Nevertheless, for the reader's convenience, we now provide a self-contained proof of this fact:

\

\begin{lemma}[Heisenberg property implies the polynomial property] \label{hppp} Suppose $\{G_n\}_{n \ge 0}$ is a chain of groups which satisfies the Heisenberg algebra property with scaling $M$. Then, for every $l \in \N$:

$$\Ind^l \Res^l = f_l(\Ind \Res)$$

where $f_l$ is the polynomial given by $f_l(X) = X(X-M)(X-2M)...(X-(l-1)M)$. In particular, $\{G_n\}_{n \ge 0}$ satisfies \hyperref[*1]{Property $(*)$}.

\end{lemma} 

\begin{proof}

The statement clearly holds when $l = 1$. Suppose the statement holds for $l = k$. Then the following relation holds.
\begin{align*}
    \Ind^{k+1} \Res^{k+1} &= \Ind^k (\Ind \Res) \Res^k \\
    &= \Ind^k (\Res \Ind - M) \Res^k \\
    &= \Ind^k (\Res \Ind \Res^k - M \Res^k) \\
    &= \Ind^k (\Res (\Res \Ind - M) \Res^{k-1} - M \Res^k) \\
    &= \Ind^k (\Res^2 \Ind \Res^{k-1} - 2M \Res^k) \\
    &= \cdots \\
    &= \Ind^k (\Res^k \Ind \Res - kM \Res^k) \\
    &= \Ind^k \Res^k (\Ind \Res - kM) \\
    &= X(X - M)(X - 2M) \cdots (X - kM)
\end{align*}
\end{proof}

\

Therefore, chains of wreath products satisfy \hyperref[*1]{Property $(*)$}, and thus serve as motivating examples for the rest of the paper. For interesting computations using the power of \hyperref[*1]{Property $(*)$} when applied to wreath products, we refer the reader to Section \ref{examples}. In the next section, we will study the interplay between \hyperref[*1]{Property $(*)$} and representation theory, and the structural limitations that representation theory provides to chains satisfying \hyperref[*1]{Property $(*)$}. The first few results from the next section will be used in the computations of Section \ref{examples}.

\

\section{Constraints on surjective chains satisfying \hyperref[*1]{Property $(*)$}}

The goal of this section is to show that \hyperref[sc]{surjective chains} that satisfy \hyperref[*1]{Property $(*)$} behave very rigidly, and therefore \hyperref[*1]{Property $(*)$} is a very rare property. The chain of wreath products is therefore a very special chain. A trivial example of such a chain of groups is the constant chain at a fixed group $G$, i.e set $G_n = G$ for all $n$. Here, induction and restriction are both the identity and so clearly $\Ind^l \Res^l = \Ind \Res$. Thus, if we set $f_l(X) = X$ for all $l$, then we have that $\Ind^l \Res^l = f_l(\Ind \Res)$. However, this example is not very interesting, and we will henceforth restrict our attention to non-constant chains which satisfy \hyperref[*1]{Property $(*)$}. 

\

\subsection{Character theoretic properties of surjective chains of groups}

\

We will now present a theorem that gives a formula for character table columns in terms of the representation theory of smaller groups of our chain and the $\Ind \Res$ operator. It is incredibly useful from a computational perspective, but moreover, it will serve as a foundational result for the rest of the paper. Though the idea involves original ideas, the proof is short and self-contained. 

\

\begin{theorem} \label{ecifpp}
\normalfont{(Extracting character information from the polynomial property)}. Let $\{G_n\}_{n \ge 0}$ be a \hyperref[sc]{surjective chain} of groups and $l \in \N$ is such that $\Ind^l \Res^l = f_l(\Ind \Res)$ for some polynomial $f_l$. Let $\alpha \in G_k \inj G_{k+l}$, and let $\delta_\alpha \in R(G_{k+l})$ be the character column of $\alpha$, that is: $\delta_\alpha = \sum_{\textrm{irreps $u$ of $G_{k+l}$}} \chi_u(\alpha) u$. For each representation $w$ of $G_k$, let $\tilde{w}$ denote a lift to $G_{k+l}$, i.e an element of $R(G_{k+l})$ such that $\Res^l(\tilde{w}) = w$. Then:

$$\delta_\alpha = f_l(\Ind \Res) \left(\sum_{\textrm{irreps $w$ of $G_k$}} \chi_w(\alpha) \cdot \tilde{w} \right)$$
\end{theorem}

\

\begin{proof}
Observe that since $\alpha \in G_k \inj G_n$, and since restrictions preserve character values, we have that for every $u \in R(G_n)$, $\chi_u(\alpha) = \chi_{\Res u}(\alpha)$. Now we replace this global restriction operator with the local $\Ind \Res$ operator as follows: 

\

\begin{equation*}
\begin{split}
\chi_{\Res^{n-k} u}(\alpha) &= \sum_{\textrm{ irreps $w$ of $G_k$}} \inner{ \Res^{n-k} u}{w} \chi_w(\alpha) \\ 
&= \sum_{\textrm{ irreps $w$ of $G_k$}} \chi_w(\alpha)  \inner{u}{\Ind^{n-k} w} \\ 
&= \sum_{\textrm{ irreps $w$ of $G_k$}} \chi_w(\alpha) \inner{u}{\Ind^{n-k} \Res^{n-k} \tilde{w} } \\ 
&= \sum_{\textrm{ irreps $w$ of $G_k$}} \chi_w(\alpha) \inner{u}{f_l(\Ind \Res) \tilde{w} } \\ 
&= \inner{f_l(\Ind \Res) \left( \sum_{\textrm{ irreps $w$ of $G_k$}} \chi_w(\alpha) \tilde{w} \right) }{u} \\ 
\end{split}
\end{equation*}

\

In the above, the second line followed from Frobenius reciprocity, the third line was from the definition of the lift $\tilde{w}$, and the fourth line was from our hypothesis.

\

We therefore have that $\inner{\delta_\alpha}{u} = \inner{\sum_{\textrm{ irreducible representations $w$ of $G_k$}} \chi_w(\alpha) f_l(\Ind \Res) \tilde{w} }{u}$ for all $u$. We therefore conclude that: 

\

$$\delta_\alpha = f_l(\Ind \Res) \left(\sum_{\textrm{ irreducible representations $w$ of $G_k$}} \chi_w(\alpha) \tilde{w} \right)$$

which proves the result.
\end{proof}

\

\begin{adjustwidth}{1cm}{}

\textbf{A couple of remarks about the theorem:} 

\

(1) \underline{Deducing character theory via the above theorem:} The theorem reduces the problem of computing the character column $\delta_\alpha$ of $G_n$ to a problem in $G_k$. The $\Ind \Res$ operator allows us to linear algebraically compute character columns via the provided formulae. We see that information about $G_k$ appears in the formula via the lifts $\tilde{w}$ - we only need to lift representations of $G_k$ (which are a relatively small number of representations for $k << n$). There are {systematic lifting procedures} that enable one to easily construct lifts $\tilde{w}$ for symmetric groups and wreath products, see Section \ref{examples} for further discussions.
\

(2) \underline{This theorem can be thought of as a way to compute characters using ring structure:} Because $U \x Ind(W) = \Ind(\Res(U) \x W)$, we can take $W = t$ (the trivial representation of the subgroup) and $U = u$ (some arbitrary representation of our group) to obtain that: $u \x \Ind(t) = \Ind \Res(u)$. That is:
\end{adjustwidth}

\begin{adjustwidth}{1cm}{}
\begin{observation}[] \label{tobs} 
Let $H \le G$ be a subgroup and let $t$ be the trivial representation of $H$. Then we have the following equality of operators: $$\Ind_H^G \Res_H^G = \Ind_H^G(t) \x$$
\end{observation} 

Therefore, in light of the above observation, $\Ind \Res$ can be thought of as the element $\Ind(t)$ in the representation ring $(R(G_n), \x)$. Our theorem gives expressions for the character table columns in terms of a polynomial in $\Ind \Res$ which translates to sums of tensor powers of $\Ind(t)$. 
\end{adjustwidth}

\

As one may observe, the key component of Theorem \ref{ecifpp} is the idea of taking lifts of representations. Given an irreducible representation $w$ in $R(G_k)$, we want to construct a lift $\tilde{w} \in R(G_n)$ such that $\Res^{n-k} \tilde{w} = w$. This is possible because, for \hyperref[sc]{surjective chains}, restriction is surjective by definition. Note that symmetric groups and wreath products form surjective chains since the branching rule is surjective. For such chains, it is possible to lift an irreducible representation from $R(G_k)$ to $R(G_n)$ for all $n$. This is not a cumbersome process since there are only $k!$ representations to lift, which is relatively insignificant when $k << n$.

\

\subsection{Proving the conjugacy class constraint} \label{pccc}

\

The goal of this subsection is to prove a family of constraints that our polynomials $\{f_l\}_{l \in \N}$ must satisfy that is indexed by conjugacy classes. We will call the constraint corresponding to a conjugacy class $[h]$ as the $[h]$-class constraint. One constraint that will turn out to be extremely important to us is the $[e]$-class constraint, which we will simply call ``the identity class constraint''. The feature of our method which allows us to recover conjugacy class information is that our method allows us to directly obtain the global structure of the character column, rather than obtaining individual character values separately. We begin with a Lemma which states that each group in the chain controls its own fusion, i.e. elements of a subgroup are conjugate in the subgroup if and only if they are conjugate in the group:

\

\begin{lemma} \label{ccclemma}
Let $\{G_n\}_{n \ge 0}$ be a surjective chain that satisfies \hyperref[*1]{Property $(*)$}. Let $h \in G_{n-l}$ for some $l < n \in \N$. Then, $f_l$ satisfies the conjugacy class constraint:

$$f_l \left(\frac{|G_n| \cdot |[h]_n \intersection G_{n-1}|}{|G_{n-1}| \cdot |[h]_n|}\right) = f_l(\chi_{\Ind(t)}(h)) =  \frac{|G_n| \cdot |[h]_{n-l}|}{|G_{n-l}|  \cdot |[h]_n|} $$
\end{lemma}

\begin{proof} For ease of notation, set $k = n-l$. By Theorem \ref{ecifpp}, 

\

$$\delta_h = f_l(\Ind \Res) \left(\sum_{\textrm{irreducible representations $w$ of $G_k$}} \chi_w(h) \tilde{w}\right)$$

\

Write $X = \Ind(t) \x = \Ind \Res$ (by Observation \ref{tobs}). Since the trivial representation $t$ has real character, $\Ind(t)$ also has real character and is therefore self-dual. The adjoint of a representation with respect to the $\Hom( \cdot, \cdot)$ bilinear form and the $\x$ multiplication is the dual representation. Therefore, by self-duality of $X$, we know that $\inner{Xu_1}{ u_2} = \inner{u_1}{ X u_2}$ for any $u_1, u_2 \in R(G_n)$. Knowing that conjugacy classes can be recovered from $L^2$ norms, we compute: 

\

\begin{equation*}
\begin{split}
||\delta_{h}||^2 
&= \inner{f_l(X) \left( \sum_{\textrm{irreps $w$ of $G_k$}} \chi_w(h) \tilde{w}\right)}{\delta_{h}} \\ 
&= \sum_{\textrm{irreps $w$ of $G_k$}} \chi_w(h) \inner{f_l(X) \tilde{w}}{\delta_{h}} \\ 
&= \sum_{\textrm{irreps $w$ of $G_k$}} \chi_w(h) \inner{\tilde{w}}{f_l(X) \delta_{h}} \\ 
&= \sum_{\textrm{irreps $w$ of $G_k$}} \chi_w(h) f_l(\chi_X(h)) \bar{\chi_{\tilde{w}}(h)} \\ 
&= \sum_{\textrm{irreps $w$ of $G_k$}} \chi_w(h) f_l(\chi_X(h)) \bar{\chi_w(h)} \\ 
&= f_l(\chi_X(h)) \cdot \sum_{\textrm{irreps $w$ of $G_k$}} |\chi_w(h)|^2 \\
&= f_l(\chi_X(h)) \cdot ||\delta^{(k)}_{h}||^2 \\
\end{split}    
\end{equation*}

where $\delta^{(k)}_{h}$ is the character column of $G_k$ corresponding to $h$, i.e $\delta^{(k)}_{h} = \sum_{\textrm{irreps $w$ of $G_k$}} \chi_w(h) w \in R(G_k)$. In the above computation, the first line followed from Theorem \ref{ecifpp}, the third line followed from self-duality of $X$, and the fifth line followed from the fact that $\Res(\tilde{w}) = w$ and so their characters coincide at $h$.

\

Translating $L^2$ norms into conjugacy classes, we see that: 

\begin{equation*}
\begin{split}
f_l(\chi_X(h)) &= \frac{||\delta_{h}||^2}{||\delta^{(k)}_{h}||^2} \\
&= \frac{\left(\frac{|G_n|}{|[h]_n|}\right)}{\left(\frac{|G_k|}{|[h]_k|}\right)} \\
&= \frac{|G_n| \cdot |[h]_k|}{|G_k|  \cdot |[h]_n|} \\
\end{split}    
\end{equation*}

\

To complete the proof, simply note that since $X = \Ind(t)$, the induced character formula tells us that:

\begin{equation*}
\begin{split}
\chi_{\Ind(t)}(h) &= \frac{1}{|G_{n-1}|} \sum_{g \in G_n \textrm{ s.t } ghg^{-1} \in G_{n-1}} \chi_t(ghg^{-1}) \\
&= \frac{1}{|G_{n-1}|} \sum_{h' \in [h]_n \intersection G_{n-1}} |\{g \in G_n \sep g h' g^{-1} = h' \}| \cdot 1 \\
&= \frac{1}{|G_{n-1}|} \sum_{h' \in [h]_n \intersection G_{n-1}} |C_{G_n}(h')| \\
&= \frac{1}{|G_{n-1}|} \sum_{h' \in [h]_n \intersection G_{n-1}} \frac{|G_n|}{| [h']_{n}|} \textrm{ (by orbit-stabilizer) } \\
&= \frac{|G_n|}{|G_{n-1}|} \cdot \frac{|[h]_n \intersection G_{n-1}|}{|[h]_{n}|} \\
\end{split}    
\end{equation*}

\

This completes the proof of the Lemma.
\end{proof}

\

We now prove a theorem which is a slight change from the Lemma above. We wish to replace $[h]_n \intersection G_{n-1}$ with $[h]_{n-1}$. That is, we require $G_{n-1}$ to control its own fusion in $G_n$.

\begin{theorem} \label{ccc}
\normalfont{(The conjugacy class constraint).} Let $\{G_n\}_{n \ge 0}$ be a non-constant \hyperref[sc]{surjective chain} of groups which satisfies \hyperref[*1]{Property $(*)$}; that is, there exists polynomials $\{f_l\}_{l \ge 1}$ such that $\Ind^l \Res^l = f_l(\Ind \Res)$. Then, for any $n \ge l$ and $h \in G_{n-l}$, $f_l$ satisfies the conjugacy class constraint

$$f_l \left(\frac{|G_n| \cdot |[h]_{n-1}|}{|G_{n-1}| \cdot |[h]_n|}\right) = f_l(\chi_{\Ind(t)}(h)) =  \frac{|G_n| \cdot |[h]_{n-l}|}{|G_{n-l}|  \cdot |[h]_n|} $$

where $[h]_j$ denotes the conjugacy class of $h \in G_j$.
\end{theorem}

\begin{proof}
By Lemma \ref{ccclemma}, it suffices to prove that $[h]_n \intersection G_{n-1} = [h]_{n-1}$. Note that if two elements of $G_{n-1}$ are conjugate in $G_{n-1}$, they are automatically conjugate in $G_{n}$ and so $[h]_{n-1} \subseteq [h]_n \intersection G_{n-1}$. We need to show that this is actually an equality. Suppose it were not. Then, there is some $h' \in G_{n-1}$ such that $h' \in [h]_n \intersection G_{n-1}$ but $h' \not\in [h]_{n-1}$. Since $h' \not\in [h]_{n-1}$, choose any class function $\phi: G_{n-1} \to \C$ such that $\phi(h) \not= \phi(h')$. This class function cannot be in the image of $\Res$, since for any class function $\psi: G_n \to \C$, $\Res(\psi)(h') =\psi(h') = \psi(h) = \Res(\psi)(h)$. Thus, $\Res$ is not surjective, which is a contradiction.
\end{proof}

\

Out of all these constraints, the most important constraint is the identity class constraint, obtained by taking $h = e$, and it is worth stating it separately: 

\

\begin{proposition}[The identity class constraint] \label{idcc}
Let $\{G_n\}_{n \ge 0}$ be a surjective chain of groups that satisfies \hyperref[*1]{Property $(*)$}. For each $n \in \N$, let $a_n = \frac{|G_n|}{|G_{n-1}|}$. Then:
 
$$f_l(a_n) = a_n a_{n-1} ... a_{n-l+1}$$ 
\end{proposition}

\begin{proof}
Setting $h = e$ in Theorem \ref{ccc}, and observing that $[e]_m = 1$ for all $m$, we see that:

\

$$f_l \left( \frac{|G_n|}{|G_{n-1}|} \right) = \frac{|G_n|}{|G_{n-l}|}$$

\

which can be rewritten as:

$$f_l(a_n) = a_n a_{n-1} ... a_{n-l+1}$$ 
\end{proof}

\

\underline{Recovering conjugacy class information for symmetric groups}

\

For convenience, we use $X = Ind(t)$ as the permutation representation of $G_n$ over the coset $[G_n : G_{n-1}]$.

Consider a permutation $\tau$. Let $S_k$ be the smallest symmetric group that $\tau$ lives in, i.e. the non-fixed elements of $\tau$ are relabeled as $\{1, \cdots, k\}$. The number of fixed points that $\tau$ has when acting on $\{1,...,n\}$ is $n-k$, and so $\chi_X(\tau) = n-k$. Since $(n-k)_{n-k} = (n-k)!$, and since $\frac{\# S_n}{\# S_k} = \frac{n!}{k!}$, Theorem \ref{ccc} informs us that: 

$$\#[\tau]_n = \frac{n!}{k!(n-k)!} \#[\tau]_k = \choose{n}{k} \#[\tau]_k$$

\

That is, we have recovered something that we already knew: Since the conjugates of $\tau$ are obtained by first choosing $n-k$ fixed points and then choosing an element of the cycle type of $\tau$ on the remaining $k$ points, we have that the size of the conjugacy class of $\tau$ should be the number of ways to choose $n-k$ fixed points multiplied by $\#[\tau]_k$.

\

\subsection{$\{a_n\}_{n \in \N}$ is necessarily infinite}

\

Let $\{G_n\}_{n \ge 0}$ be a non-constant \hyperref[sc]{surjective chain} that satisfies \hyperref[*1]{Property $(*)$}. Recall that $a_n := \frac{|G_n|}{|G_{n-1}|}$. Why might we want $\{a_n\}_{n \in \N}$ to be an infinite set? The reason is that the infinitude of $\{a_n\}_{n \in \N}$ allows us to use analytic arguments in conjunction with the identity class constraint to deduce very rigid properties that our polynomials must satisfy. Notice that if  $\{G_n\}_{n \ge 0}$ is a constant chain (i.e $G_n = G$ for all $n$), then $a_n = 1$ for all $n$ and $\{G_n\}_{n \ge 0}$ satisfies \hyperref[*1]{Property $(*)$} with $f_l(X) = X$ for all $l$. Therefore, the non-constancy of $\{G_n\}_{n \ge 0}$ is a prerequisite for infinitude of $\{a_n\}_{n \in \N}$. Our goal of this subsection is to prove the following Lemma:

\

\begin{lemma}[Infinitude of $\{a_n\}_{n \in \N}$] \label{infan} Let $\{G_n\}_{n \ge 0}$ be a non-constant \hyperref[sc]{surjective chain} of groups satisfying \hyperref[*1]{Property $(*)$}. Define $a_n = \frac{|G_n|}{|G_{n-1}|}$. Then: $\{a_n\}_{n \in \N}$ is an infinite set.
\end{lemma} 

\

We will prove our goal in 3 steps. 

\

\begin{lemma}[Step 1] \label{step1}
Let $\{G_n\}_{n \ge 0}$ be a non-constant surjective chain of groups such that $\Ind^2 \Res^2 = f_2(\Ind \Res)$ for some polynomial $f_2$. Suppose that $h_1, h_2 \in G_{n-1} \inj G_n$ are such that $\chi_{\Ind(t_{n-2})}(h_1) \not= \chi_{\Ind(t_{n-2})}(h_2)$. Let $t_i$ denote the trivial representation of $G_i$.Then, 

$$\chi_{\Ind(t_{n-1})}(h_1) \not= \chi_{\Ind(t_{n-1})}(h_2)$$
\end{lemma} 

\begin{proof}
By the proof of Theorem \ref{ccc}, we know that every group $G_k$ of our chain controls its own fusion. Since $f_2(X) = \Ind^2(t)$, we have that for any  $h \in G_{n-1} \inj G_n$, $f_2(\chi_{\Ind(t_{n-1})}(h)) = \chi_{\Ind^2(t_{n-2})}(h)$. By the induced character formula, we have that:

\begin{equation*}
\begin{split}
f_2(\chi_{\Ind(t_{n-1})}(h)) &= \chi_{\Ind^2(t_{n-2})}(h) \\ 
&= \chi_{\Ind(\Ind(t_{n-2}))}(h) \\ 
&= \frac{1}{|G_{n-1}|} \sum_{g \in G_n \textrm{ such that } g h g^{-1} \in G_{n-1}} \chi_{\Ind(t_{n-2})} (g h g^{-1}) \\ 
&= \frac{1}{|G_{n-1}|} \sum_{g \in G_n \textrm{ such that } g h g^{-1} \in G_{n-1}} \chi_{\Ind(t_{n-2})} (h) \\ 
&= \chi_{\Ind(t_{n-2})} (h) \cdot \frac{1}{|G_{n-1}|} \sum_{g \in G_n \textrm{ such that } g h g^{-1} \in G_{n-1}} 1 \\ 
&= \chi_{\Ind(t_{n-2})} (h) \cdot \chi_{\Ind(t_{n-1})}(h) \\ 
\end{split}
\end{equation*}

\

In the above computation, the third line followed from the induced character formula, and the fourth line followed from the fact that $G_{n-1}$ controls its own fusion. 

\

Now, suppose $h_1, h_2 \in G_{n-1} \inj G_n$ such that  $\chi_{\Ind(t_{n-1})}(h_1) = \chi_{\Ind(t_{n-1})}(h_2)$. Observe that $\chi_{\Ind(t_{n-1})}(h_1) \neq 0$ because $h_1$ fixes $eG_{n-1} \in G_n/G_{n-1}$. Then, $\frac{f_2(\chi_{\Ind(t_{n-1})}(h_1))}{\chi_{\Ind(t_{n-1})}(h_1)} = \frac{f_2(\chi_{\Ind(t_{n-1})}(h_2))}{\chi_{\Ind(t_{n-1})}(h_2)}$. So, we conclude that $\chi_{\Ind(t_{n-2})}(h_1) =  \chi_{\Ind(t_{n-2})}(h_2)$, which proves the contrapositive statement of the Lemma and therefore completes the proof. 
\end{proof}

\

\begin{lemma}[Step 2] \label{step2} Let $\{G_n\}_{n \ge 0}$ be a non-constant surjective chain of groups such that $\Ind^2 \Res^2 = f_2(\Ind \Res)$ for some polynomial $f_2$. Then, $\Ind(t) \in R(G_n)$ has at least $n$ distinct character values. 
\end{lemma} 

\begin{proof}
We proceed by induction on $n$. For $n = 1$ there is nothing to show since every representation must have at least one character value. 

\

Suppose it is true that $\Ind(t_{n-2})$ has $n-1$ distinct character values, say $\chi_{\Ind(t_{n-2})}(h_1), \chi_{\Ind(t_{n-2})}(h_2), ... , \chi_{\Ind(t_{n-2})}(h_{n-1})$. By \hyperref[step1]{Lemma \ref{step1}}, $\chi_{\Ind(t_{n-1})}(h_1), \chi_{\Ind(t_{n-1})}(h_2), ... , \chi_{\Ind(t_{n-1})}(h_{n-1})$ are $n-1$ distinct character values of $\Ind(t_{n-1})$. Moreover, each of these character values are nonzero because if $h_i \in G_{n-1}$, then $h_i$ fixes the identity coset. By Burnside's Lemma, because the action of $G_n$ on $\left\{\quotient{G_n}{G_{n-1}}\right\}$ is a transitive action, there must exist some $g \in G_n$ such that $\chi_{\Ind(t_{n-1})}(g) = 0$. Thus, $0, \chi_{\Ind(t_{n-1})}(h_1), \chi_{\Ind(t_{n-1})}(h_2), ... , \chi_{\Ind(t_{n-1})}(h_{n-1})$ are $n$ distinct character values of $\Ind(t_{n-1})$ and we have therefore proved the result.
\end{proof}

\

\begin{observation}[Step 3] \label{step3}
The number of distinct character values of $\Ind(t)$ is at most $a_n$
\end{observation} 
\begin{proof}
Since the character values of $\Ind(t)$ are the number of fixed points, they are integers. Moreover $\dim(\Ind(t)) = a_n$ is the maximal character value
\end{proof}

\

Using Steps \hyperref[step1]{1}, \hyperref[step2]{2}, and \hyperref[step3]{3}, we can prove the infinitude of $\{a_n\}_{n \in \N}$ which is the main goal of this section. Recall the statement: 

\

\begin{replemma}{infan} \normalfont{(Infinitude of $\{a_n\}_{n \in \N}$).} Let $\{G_n\}_{n \ge 0}$ be a non-constant \hyperref[sc]{surjective chain} of groups satisfying \hyperref[*1]{Property $(*)$}. Define $a_n = \frac{|G_n|}{|G_{n-1}|}$. Then: $\{a_n\}_{n \in \N}$ is an infinite set.
\end{replemma} 

\begin{proof}
By \hyperref[step2]{Lemma \ref{step2}}, we have that the number of distinct character values of $\Ind(t) \in R(G_n)$ is $\ge n$. But by \hyperref[step3]{Lemma \ref{step3}}, the number of distinct character values of $\Ind(t) \in R(G_n)$ is $\le a_n$. This is only compatible if $a_n \ge n$ for every $n$. In particular, $a_n \to \infty$ as $n \to \infty$. 
\end{proof}

\

\subsection{Description of $f_l$ in terms of the two parameters $B$ and $C$}

\

Now that we have established that \hyperref[infan]{$\{a_n\}_{n \in \N}$ is an infinite collection}, we are ready to use the identity class constraint (Proposition \ref{idcc}) together with analytic methods to deduce results. The end goal of this subsection is to prove Theorem \ref{pdtp}. We begin with a simple observation:

\begin{observation}[Bounding degrees of $f_l$] \label{deg1} Let $\{G_n\}_{n \ge 0}$ be a non-constant \hyperref[sc]{surjective chain} of groups satisfying \hyperref[*1]{Property $(*)$}. Then, $f_l$ has degree at most $l$.
\end{observation} 

\begin{proof}
We first prove that $f_2$ has degree at most $2$. Suppose there are infinitely many $n \in \N$ such that $a_{n-1} > a_n$. Then $a_n^2 < f_2(a_n) = a_n a_{n-1}$
for infinitely many $n$. Because $a_n \to \infty$ as $n \to \infty$, there also exist infinitely many $n \in \N$ such that $a_{n-1} \leq a_n$. For such $n$'s, $f_2(a_n) \leq a_n^2$. Then the polynomial $f_2(X) - X^2$ has infinitely many local extrema, a contradiction. Thus there exists only finitely many $n \in \N$'s such that $a_{n-1} > a_n$. In particular, there exists a fixed positive integer $M$ such that for all $n \geq M$, $a_{n-1} \leq a_n$. Since $f_2(a_n) = a_n a_{n-1} \leq a_n^2$ for all sufficiently large $n$, we obtain that $f_2$ has degree at most $2$.

The proof for $f_l$ immediately follows because for all $n \geq M$, $f_l(a_n) = a_n \cdots a_{n-l+1} \leq a_n^l$.
\end{proof}

\

We also deduce the following Lemma, which is an interesting and important stepping stone: 

\

\begin{lemma}[Smaller polynomials divide larger polynomials] \label{divp}
Let $\{G_n\}_{n \ge 0}$ be a non-constant \hyperref[sc]{surjective chain} of groups satisfying \hyperref[*1]{Property $(*)$}. Then: $f_{l-1}$ divides $f_l$ for every $l$. 
\end{lemma} 

\begin{proof}
Recall that $a_n = \frac{|G_n|}{|G_{n-1}|}$. By the \hyperref[idcc]{Identity class constraint (Proposition \ref{idcc})}: 

$$f_l(a_n) = a_n a_{n-1} ... a_{n-l+1} = f_{l-1}(a_n)a_{n-l+1}$$

We therefore have that, for all $n \in \N$, $\frac{f_{l}(a_n)}{f_{l-1}(a_n)} = a_{n-l+1} \in \Z$. The condition that $a_{n-l+1} \in \N$ is just a consequence of Lagrange's theorem since $G_{n-l} \le G_{n-l+1}$.

\

By the divisor theorem on $\Q[X]$, we can write $f_l = f_{l-1}q + r$ where $r$ is a polynomial of degree less than $f_{l-1}$. Thus, we may write: 

\ 

$$\frac{f_l}{f_{l-1}} = q + \frac{r}{f_l}$$

\

Suppose for contradiction that $r$ is nonzero. 

\ 

Choose $\epsilon > 0$ such that for all integers $a \in \Z$, either $q(a) \in \Z$ or $|q(a) - \Z| \ge \epsilon$. Let us first show that such an $\epsilon$ exists. For pedagogical reasons, first consider the case in which $q = \frac{c}{d} X^p$ for $c,d \in \Z$. Observe that for $a \in \Z$ either $q(a) \in \Z$ or $|q(a) - \Z| \ge \frac{1}{d}$. For such a polynomial $q$, we may take $\epsilon = \frac{1}{d}$. In general, write $q = \frac{c_p}{d_p}X^p +\frac{c_{p-1}}{d_{p-1}}X^{p-1}+...+\frac{c_0}{d_0}$. For such a polynomial, we may take $\epsilon = \frac{1}{\prod_{i=0}^p d_p}$.

\ 

Choose $\epsilon$ as in the previous paragraph.  For large values of $n$, (say $n > N$), we can assert that $0 < |(\frac{r}{f_l})(a_{n})| < \epsilon$. The first inequality is simply because $r$ is nonzero. The second inequality is because the degree of $f_l$ is strictly larger than the degree of $r$, and \hyperref[infan]{$a_n \to \infty$ as $n \to \infty$ (Lemma \ref{infan})}.

\ 

Take any $n > N$. Then we have that: either  $q(a_n) \in \Z$ or $|q(a_n) - \Z| \ge \epsilon$, and $0 < (\frac{r}{f_l})(a_n) < \epsilon$. Therefore, we have that:

$$(q + \frac{r}{f_l})(a_n) \not\in \Z$$

Therefore, $$(\frac{f_l}{f_{l-1}})(a_n) \not\in \Z$$ for this $n$, and this is a contradiction since we remarked earlier that $\frac{f_l}{f_{l-1}}(a_n) \in \Z$ for all $n$. We therefore must have that $r = 0$, i.e $f_{l-1}$ divides $f_l$.
\end{proof}

\

An immediate and useful application of the two above results is the following:

\

\begin{observation}[Determining $f_2$] \label{f2}
Let $\{G_n\}_{n \ge 0}$ be a non-constant \hyperref[sc]{surjective chain} of groups satisfying \hyperref[*1]{Property $(*)$}. Then, there exist constants $A, C \in \Q$ such that $$f_2(X) = AX(X-C)$$
\end{observation} 

\begin{proof}
Since $\Ind^1 \Res^1 = \Ind \Res$, we have that $f_1(X) = X$. By the Lemma \ref{divp}, $f_1$ must divide $f_2$, and by Observation \ref{deg1}, $f_2$ is either a degree $1$  polynomial or a degree $2$ polynomial. Thus, we must have that either $f_2(X) = AX$, or $f_2(X) = AX(X-C)$ for some constants $A, C \in \Q$. If the former is true, then by the identity class constraint (Proposition \ref{idcc}), $a_n a_{n-1} = f_2(a_n) = A a_n$ for all $n$, and thus, $a_{n-1} = A$ is a constant for all $n$, thereby contradicting the infinitude of $\{a_n\}_{n \in \N}$  (Lemma \ref{infan}). Thus, $$f_2(X) = AX(X-C)$$. 
\end{proof}

\

We therefore observe the following recursion formula:

\

\begin{observation}[Recursion of $\{a_n\}_{n \in \N}$] \label{recan}
Let $\{G_n\}_{n \ge 0}$ be a non-constant \hyperref[sc]{surjective chain} of groups satisfying \hyperref[*1]{Property $(*)$}. Then, there exist constants $A, C \in \Q$ such that $$a_n = \frac{1}{A}a_{n-1} + C$$
\end{observation} 

\begin{proof}
Using the identity lass constraint (Proposition \ref{idcc}) for $f_2$ above, we see that $a_n a_{n-1} = A a_n (a_n - C)$. The observation follows.
\end{proof}

\

Thus far, we do not have any restrictions on $A,C \in \Q$. The following proposition shows that $\frac{1}{A}$ and $C$ must both be integers.  

\

\begin{proposition}[$C, \frac{1}{A} \in \Z$] \label{cbint}
Suppose $\{a_n\}$ is a sequence of positive integers such that $a_n = \frac{1}{A} a_{n-1} + C$ for fixed rational numbers $A$ and $C$. Suppose also that $a_n \to \infty$ as $n \to \infty$. Then $A = \frac{1}{B}$ for some positive integer $B$ and an integer $C$.
\end{proposition} 

\begin{proof}
Write $A = \frac{p}{q}$ where $(p,q) = 1$, and write $C = \frac{c}{d}$ for $(c,d) = 1$. By repeatedly applying the identity $a_n = \frac{1}{A} a_{n-1} + C$, we observe that $a_n$ can be expressed as the following equation in terms of $a_1$: 

\begin{align*}
    a_n &= \frac{1}{A^{n-1}} \left( a_1 + CA + CA^2 + \cdots + CA^{n-1} \right)\\
    &= \frac{q^{n-1}}{p^{n-1}} \left( a_1 + C \left( \frac{p}{q} + \cdots + \frac{p^{n-1}}{q^{n-1}} \right) \right)
\end{align*}

\

We then obtain that: 

\begin{equation*}
    a_n d p^{n-1} = a_1 d q^{n-1} + c p q^{n-2} + \cdots + cp^{n-1}
\end{equation*}

Now, suppose that $p \geq 2$. Since $(p,q) = 1$, taking the above equation mod $p^i$ for all $0 \le i \le n-1$ and dividing by suitable powers of $q^j$ gives:

\begin{align*}
    a_1 d &\equiv 0 \; \text{mod} \; p \\
    a_1 d q &\equiv -cp \; \text{mod} \; p^2 \\
    a_1 d q^2 &\equiv -c (p^2 + pq) \; \text{mod} \; p^3 \\
    &\vdots \\
    a_1 d q^{n-2} &\equiv -c (p^{n-2} + \cdots + pq) \; \text{mod} \; p^{n-1}
\end{align*}

Since $q^j$ is a unit in $\Z/p^i\Z$, we also obtain: 

\begin{align*}
    a_1 d &\equiv 0 \; \text{mod} \; p \\
    a_1 d &\equiv -c \frac{p}{q} \; \text{mod} \; p^2 \\
    a_1 d &\equiv -c \left( \frac{p^2}{q^2} + \frac{p}{q} \right) \; \text{mod} \; p^3 \\
    &\vdots \\
    a_1 d &\equiv -c \left( \frac{p^{n-2}}{q^{n-2}} + \cdots + \frac{p}{q} \right) \; \text{mod} \; p^{n-1}
\end{align*}

But notice that $a_1 d$ is a fixed integer. Then there exists an integer $k$ such that $a_1 d< p^k$. Therefore, taking sufficiently large $n$, we obtain that 

\begin{equation*}
    -c \left( \frac{p^{k-2}}{q^{k-2}} + \cdots + \frac{p}{q} \right) \equiv -c \left( \frac{p^{k-1}}{q^{k-1}} + \frac{p^{k-2}}{q^{k-2}} + \cdots + \frac{p}{q} \right) \; \text{mod} \; p^n
\end{equation*}

\

In particular, 

\

\begin{equation*}
    c \frac{p^{k-1}}{q^{k-1}} \equiv 0 \; \text{mod} \; p^n
\end{equation*}

\

Hence, $c \equiv 0 \; \text{mod} \; p^{n-k-1}$. 

\

From this relation, we first show that $A = \frac{1}{q}$. If $c = 0$, then $C = 0$. Therefore, $a_n = \frac{1}{A} a_{n-1} = \frac{q}{p} a_{n-1}$. This is a contradiction because $a_n$ is a positive integer for all $n \in \N$. If $c \neq 0$. Then the equivalence $c \equiv 0 \; \text{mod} \; p^{n-k-1}$ is a contradiction because it holds for arbitrarily large $n$ while $c$ is a fixed integer. 

\

Hence, $p = 1$ so $A = \frac{1}{q}$ for some integer $q \in \mathbb{Z}$. We can fix $B = q$. Note that we in fact require that $B$ is positive because $\{a_n\}_{n \in \N}$ is a sequence of positive integers. Next, we recall that $a_n = \frac{1}{A}a_{n-1} + C$. Since $a_n, a_{n-1}$, and $\frac{1}{A} = B$ are integers, we require that $C$ is an integer.
\end{proof}

\

If we concatenate the theory we have built thus far, we notice that we have proved Part (1) of Theorem \ref{pdtp}. We may use this to strengthen Observation \ref{deg1} as follows:

\

\begin{lemma}[$f_l$ has degree $l$] \label{deg2}
Let $\{G_n\}_{n \ge 0}$ be a non-constant \hyperref[sc]{surjective chain} of groups satisfying \hyperref[*1]{Property $(*)$}. Then, $f_l$ is a polynomial of degree $l$.
\end{lemma} 

\begin{proof}
By the \hyperref[idcc]{Identity class constraint (Proposition \ref{idcc})}, 

$$f_l(a_n) = a_n a_{n-1} ... a_{n-l+1} = a_n \cdot a_n \frac{a_{n-1}}{a_n} \cdot  a_n \frac{a_{n-2}}{a_{n-1}} \frac{a_{n-1}}{a_{n}} ... \cdot a_n \frac{a_{n-l+1}}{a_{n-l+2}} \frac{a_{n-l+2}}{a_{n-l+3}} .. \frac{a_{n-1}}{a_{n}}$$

Now observe that by Part (1) of Theorem \ref{pdtp}, $\frac{a_m}{a_{m+1}} = \frac{a_m}{Ba_m+C} = \frac{1}{B + \frac{C}{a_m}}$ for every $m$, and thus we obtain that for $m$ large:

$$\frac{1}{B}-\epsilon \le \frac{a_m}{a_{m+1}} \le \frac{1}{B} + \epsilon$$

In particular, for $m$ large:

$$\left( \frac{1}{B}-\epsilon \right)^{\frac{(l-1)(l)}{2}} a_n^l \le f_l(a_n) \le \left( \frac{1}{B}+\epsilon \right)^{\frac{(l-1)(l)}{2}} a_n^l$$

and thus $f_l$ has degree $l$ by the \hyperref[infan]{Infinitude of $\{a_n\}_{n \in \N}$ (Lemma \ref{infan})}.
\end{proof}

\

We are now ready to prove Theorem \ref{pdtp} in its full entirety. For the reader's convenience, we will state the theorem below:

\

\begin{reptheorem}{pdtp} [Polynomials are determined by two parameters] 
Let $\{G_n\}_{n \ge 0}$ be a non-constant \hyperref[sc]{surjective chain} of groups which satisfies \hyperref[*1]{Property $(*)$}; that is, there exists polynomials $\{f_l\}_{l \ge 1}$ such that $\Ind^l \Res^l = f_l(\Ind \Res)$.  Let $a_n = \frac{|G_n|}{|G_{n-1}|}$. Then: \\
\begin{addmargin}[2em]{0em}
(1) $\{a_n\}_{n \in \N}$ is an infinite set. Moreover, there exists two parameters $B \in \N$ and $C \in \Z$ such that for all $n$, $a_n = B a_{n-1} + C$. \\
(2) The two parameters $B$ and $C$ from Part (2) completely determine the polynomials $\{f_l\}_{l \in \N}$. $f_l$ can be expressed in terms of these parameters as: $$f_l(X) = \frac{1}{B^{\frac{l(l-1)}{2}}} X (X - C) (X - [1+ B]C) (X - [1 + B + B^2]C) ... (X - [1 + B +  .... +B^{l-2} ] C) $$
\end{addmargin}
\end{reptheorem}

\begin{proof}
(1) follows directly from Lemma \ref{infan}, Observation \ref{recan} and Proposition \ref{cbint}. It therefore remains to prove (2). 

\

We prove the formula in (2) by induction on $l$. When $l = 1$, we have $f_1(X) = X$, which satisfies the formula prescribed by the proposition. Suppose that $f_l$ satisfies the formula prescribed by the proposition. Since \hyperref[divp]{$f_l$ divides $f_{l+1}$ (Lemma \ref{divp})} and \hyperref[deg2]{$f_l$ has degree $l$ and $f_{l+1}$ has degree $l+1$ (Lemma \ref{deg2})}, we may write: 

$$\frac{f_{l+1}(X)}{f_l(X)} = B_{l+1} (X - C_{l+1})$$

for some $B_{l+1}, C_{l+1} \in \Q$. 

\

Evaluating the above expression at $a_n$ and applying the \hyperref[idcc]{Identity class constraint (Proposition \ref{idcc})} to simplify the left hand side, we see that:  
 
\begin{equation*}
    a_{n-l} = B_{l+1} (a_n - C_{l+1})
\end{equation*}

Since $a_{n-1} = \frac{1}{B} (a_{n} - C)$, we have that: 

\begin{equation*}
    a_{n-l} = \frac{1}{B^l} \left(a_n - C(1 + B + \cdots + B^{l-1}) \right)
\end{equation*}

\

Comparing expressions, we see that $B_{l+1} = \frac{1}{B^l}$ and $C_{l+1} = C(1 + B + \cdots + B^{l-1})$. Thus: 

\begin{equation*}
\begin{split}
f_{l+1}(X) &= \frac{1}{B^l} \cdot \frac{1}{B^{\frac{l(l-1)}{2}}} X (X - C) (X - [1+ B]C) (X - [1 + B + B^2]C) ...  (X - [1 + B +  .... +B^{l-1} ] C) \\
&= \frac{1}{B^{\frac{(l+1)(l)}{2}}} X (X - C) (X - [1+ B]C) (X - [1 + B + B^2]C) ... (X - [1 + B +  .... +B^{l-1} ] C) \\ 
\end{split}
\end{equation*}

which completes the inductive step.
\end{proof}

\

\begin{adjustwidth}{1cm}{}
\textbf{Some remarks about the theorem:} 

\

(1) \underline{Relation to results via dual graded graphs:} We note that one can prove analogous formulae by viewing Grothendieck rings of towers of algebras $\bigoplus_{n \geq 0} \mathbb{C} G_n$ as a dual graded graph, see in particular Theorem 1.1 and Theorem 7.2 of \cite{dg2} and Theorem 3.7 of \cite{dg3}. However, our result is fundamentally different, since those results assume that we have an injection $\mathbb{C}[G_k] \otimes \mathbb{C}[G_l] \to \mathbb{C}[G_n]$, which is a strong condition that automatically implies a linear recurrence relation, see Theorem 3.2 and Theorem 7.4 of \cite{dg2}. Those papers then deduce the polynomial relation from the linear recurrence relation using Section 1.4 of \cite{dg1}. In contrast, the goal of our theorem is to assume a polynomial recurrence relation, and from this deduce that a linear recurrence relation must hold, thereby addressing a converse direction. The arguments we use also have a very different flavor from those used in \cite{dg2}, \cite{dg3}, and \cite{dg1}, since our arguments are fundamentally representation theoretic, whereas those papers extensively use combinatorial properties of towers of algebras.

\ 

(2) \underline{The two parameter family:}  A consequence of Theorem \ref{pdtp} is that whenever a surjective chain satisfies \hyperref[*1]{Property $(*)$}, the polynomials $f_l$ belong to a two-parameter family of polynomials and the parameters are integers ($B$ and $C$). This is certainly a very rigid constraint. Moreover, by Lemma \ref{hppp}, we know that there exist chains for $B = 1$ and any positive integer $C$. Furthermore, any given $f_k$ determines $f_l$ because $f_k$ uniquely determines the two parameters $B$ and $C$ as shown in Part (2) of Theorem \ref{pdtp}. It is a natural question to ask whether it is possible to have $C$ to be non-positive. In fact, if $G_0$ is trivial, then $C$ must be positive. 

\

(3) \underline{Reading off polynomials from three successive values: $a_n, a_{n+1}, a_{n+2}$:} Observe that the orders of four consecutive groups $G_{n-1}, G_n, G_{n+1}, G_{n+2}$ determines $f_l$ for all $l$. This suggests that the sequence $\{f_l\}_{l \in \N}$ behaves very rigidly. Indeed, note that $a_{n+1} = B a_n + C \implies C = a_{n+1} - B a_n$, and thus, $a_{n+2} = Ba_{n+1} + C = (B+1) a_{n+1} - B a_n \implies B(a_{n+1} - a_n) = (a_{n+2} - a_{n+1})$, and thus $B$ is uniquely determined, since one can show that $a_{n+1} \not= a_n$ for non-constant chains. And likewise $C$ is uniquely determined.

\

(4) \underline{Predicting the polynomials for wreath products from dimensions:} For the chain of wreath products $\{H^n \sdp S_n\}_{n \in \N}$, $a_n = n|H|$. One may observe that $a_n = a_{n-1} + |H|$. The constants $B$ and $C$ in Theorem \ref{pdtp} are then $1$ and $|H|$ respectively, and Theorem \ref{pdtp} therefore implies that if $\{H^n \sdp S_n\}_{n \in \N}$ satisfies \hyperref[*1]{Property $(*)$}, then $f_l(X) = X(X-|H|)(X-2|H|)...(X-(l-1)|H|)$.
\end{adjustwidth}

\

As a corollary of the above result, we obtain a result which has importance in the world of dual graded graphs (see \cite{dg1} for a rigorous treatment of the subject).

\

\begin{corollary}
Let $\{G_n\}_{n \ge 0}$ be a non-constant \hyperref[sc]{surjective chain} of groups which satisfies \hyperref[*1]{Property $(*)$}. Then there exist constants $B$ and $C$ such that

$$\Res \Ind - B \Ind \Res = C \textrm{Id}$$

That is, the language of dual graded graphs (see  \cite{dg1}), $\Ind$ and $\Res$ operators have recurrent commutation relation if and only if they have linear commutation relation. 
\end{corollary}

\

\subsection{If $G_0 = \{e\}$, the roots of $f_l$ correspond to characters values of $\Ind(t) \in G_l$}

\

In this subsection, we will consider the special and important case when $G_0 = \{e\}$. The reason why this situation is special and insightful is that $G_0$ has only one representation, the trivial representation, which we always know how to lift. This allows us to apply our character column techniques from Theorem \ref{ecifpp} to deduce special information. Our goal is to prove:

\begin{theorem} \label{rpcv}
Let $\{G_n\}_{n \ge 0}$ be a non-constant \hyperref[sc]{surjective chain} of groups which satisfies \hyperref[*1]{Property $(*)$}. Suppose further that $G_0 = \{ e \}$. Then, the roots of $f_l$ are precisely the non-identity character values of $X = \Ind(t) \in R(G_l)$ (i.e the character values other than $\dim(\Ind(t))$).
\end{theorem}

\begin{adjustwidth}{1cm}{}
\textbf{A couple of remarks about the theorem:} 

\ 

(1) \underline{Positivity of $C$}: Theorem \ref{rpcv} shows that if $G_0$ is trivial, then the integer $C$ from Theorem \ref{pdtp} has to be positive. We do not know of any examples of chains for which $C$ is negative, and a reasonable question for the future would be to ask whether or not $C$ is forced to be positive.

\ 

(2) \underline{Predicting the polynomials $f_l$:} We may use Theorem \ref{rpcv} to predict the polynomials $f_l$. For example, consider the chain of symmetric groups $\{S_n\}_{n \in \N}$. Since $S_1 = \{e\}$, Theorem \ref{rpcv} informs us that the roots of the polynomial $f_l$ are precisely the non-identity eigenvalues of $X = \Ind(t) \in R(S_{l+1})$. These eigenvalues are the characters of the permutation representation of $S_{l+1}$, which are obtained by counting fixed points of the $S_{l+1}$ action on $\{1,2,...,l+1\}$. The possible number of fixed points a non-identity permutation can have belongs to $\{1,2,3,.., l-1\}$, and we therefore know that $f_l$ is a multiple of the falling factorial polynomial $X(X-1)...(X-(l-1))$. Dimension counting is enough to show that this multiple is $1$. 

\

(3) \underline{Generalization to wreath products:} The aforementioned argument can be generalized to the chain of wreath products $\{H^n \rtimes S_n\}_{n \in \mathbb{N}}$ as well. We thank the reviewer for pointing out this example. Suppose that $H$ acts on a set $S$. Let $K < H$ be a point stabilizer. Then the permutation representation $\Ind_{H^{n-1} \rtimes S_{n-1} \times K}^{H^n \rtimes S_n} t$ is the permutation representation of $H^n \rtimes S_n$ acting imprimitively on $H \times \{1, \cdots, n\}$. This is because the group $H^n \rtimes S_n$ acts imprimitively on the set $S \times \{1, \cdots, n\}$ with blocks $S \times \{k\}$ for $k \in \{1, \cdots, n\}$, and the point stabilizer for this action is $H^{n-1} \rtimes S_{n-1} \times K$. Because $H$ acts regularly on itself, the possible number of fixed points a non-identity permutation can have belongs to $\{|H|, 2|H|, \cdots, (n-2)|H|, n|H|\}$.
\end{adjustwidth}
\

We will now prove the statement of Theorem \ref{rpcv}.

\

\begin{lemma}[Non-identity character values are roots] \label{rpcv1}
Suppose that $\{G_n\}_{n \ge 0}$ is a surjective chain of groups that satisfies \hyperref[*1]{Property $(*)$} with the additional property that $G_0 = \{e\}$. Consider $\Ind(t) \in G_l$. Then, $\chi_{\Ind(t)}(\alpha)$ is a root of $f_l$ for every $\alpha \in G_l$ other than $\alpha = e$.  
\end{lemma} 

\

\begin{proof} Consider $e \in \{e\} = G_0$. Observe that Theorem \ref{ecifpp} informs us that:  

\ 

$$\delta_{e} = \chi_t(e)f_l(X)\tilde{e} = f_l(X)\tilde{t} = f_l(X)t$$

\ 

Consider another character column $\delta_\alpha$ such that $\alpha \neq e$. Observe that this is an eigenvector of $X$ with eigenvalue $\chi_{\Ind(t)}(\alpha)$. Indeed, this is a simple consequence of the fact that $\Ind \Res = \Ind(t) \x$. It might be worth mentioning that the fact that character columns are eigenvectors for $\Ind \Res$ with the above eigenvalues is also observed in \cite{dg6} and is given more prominence there. 

\

So, we then have that:

\

\begin{equation*}
\begin{split}
0 &= \inner{\delta_e}{\delta_\alpha} = \inner{f_l(X)t}{\delta_{\alpha}} = \inner{t}{f_l(\chi_{\Ind(t)}(\alpha))\delta_{\alpha}} = f_l(\chi_{\Ind(t)}(\alpha)) \delta_\alpha(t) = f_l(\chi_{\Ind(t)}(\alpha)) \\ 
\end{split}    
\end{equation*}
\end{proof}

\

It might be worth to remark that If we had used $\alpha = e$ with the above reasoning, then we have obtained that:

$$|G_l| = \inner{\delta_e}{\delta_e} = f_l(\frac{|G_l|}{|G_{l-1}|}) = f_l(a_l)$$

which is simply the \hyperref[idcc]{Identity class constraint (Proposition \ref{idcc})} in the special case where $G_0 = \{e\}$. 

\

Now we prove the converse:

\begin{lemma}[Roots are non-identity character values] \label{rpcv2}
Suppose that $\{G_n\}_{n \ge 0}$ is a surjective chain of groups that satisfies \hyperref[*1]{Property $(*)$} with the additional property that $G_0 = \{e\}$. Then, any root of $f_l$ is an eigenvalue of $\Ind(t) \in R(G_l)$.
\end{lemma} 

\begin{proof}
Let $\gamma$ be a root of $f_l(X)$. Then $f_l(X)$ is divisible by $X-\gamma$. By Theorem \ref{ecifpp}, $\delta_{e} = f_l(X)$, and thus, for any $\alpha \not= e$: 

\ 

$$0 = \inner{\delta_e}{\delta_
\alpha} = \inner{f_l(X)}{\delta_\alpha} = \inner{X-\gamma}{\frac{f_l(X)}{X-\gamma} \delta_\alpha}$$

\ 

Therefore, if $\frac{f_l(X)}{X-\gamma} \delta_\alpha$ is nonzero for some $\alpha$, then it is an eigenvector of $X$ with eigenvalue $\gamma$. 

\

If $\frac{f_l(X)}{X-\gamma} \delta_\alpha$ is the zero vector for all $\alpha \not= e$, then:

\begin{equation*}
    0 = \inner{\frac{f_l(X)}{X-\gamma} \delta_\alpha}{t} = \inner{\frac{f_l(X)}{X-\gamma}}{\delta_\alpha}
\end{equation*}

Hence, $\frac{f_l(X)}{X-\gamma}$ is orthogonal to all $\delta_\alpha$'s. Because $\{\delta_\alpha\}$ forms an orthogonal basis of the space $R(G_l)$, $\frac{f_l(X)}{X-\gamma}$ is a scalar multiple of $\delta_e$. This is impossible since we would have that $\frac{f_l(X)}{X-\gamma} = \lambda f_l(X)$ but the degrees on both sides of the equation are not equal. Therefore, there must exist some $\alpha$ such that $\frac{f_l(X)}{X-\gamma} \delta_\alpha$ is nonzero, and this is an eigenvector of $X$ with eigenvalue $\gamma$. 
\end{proof} 

\

We have thus proved Theorem \ref{rpcv}.

\

\section{Examples} \label{examples}

As examples, we present how Property \hyperref[*1]{$(*)$} can be utilized to excavate the representation theoretic properties of chains of symmetric groups and wreath products.

\subsection{Chains of symmetric groups}

The chain of symmetric groups is a chain of wreath products with $H = \{e\}$,. The chain of symmetric groups satisfies the Heisenberg algebra property (Property \ref{hap}) with scaling $1$. That is, the chain of symmetric groups satisfies Property \hyperref[*1]{$(*)$} with $f_l(X) = (X)_l$, where $(X)_l$ is the falling factorial polynomial. In this subsection, we will exhibit some formulae for symmetric group character columns using Theorem \ref{pdtp} and \ref{ecifpp}. 

\begin{theorem}
Let $X = \Ind \Res = \Ind(t) \x$ be the $\Ind \Res$ operator on the chain $\{S_n\}$.
\begin{enumerate}
    \item The polynomial $f_l$ is given by
    $$ f_l = X(X-1)(X-2)...(X-(l-1))$$
    \item The non-identity character values of $\Ind t$ in the symmetric group $S_n$ is the set $\{0, 1, \cdots, n-1\}$.
    \item Let $\tau \in S_k \inj S_n$. The character column $\delta_\tau$ as an element of $R(S_n)$ is given by: $$\delta_\tau = X(X-1)...(X-(n-k)+1) \left(\sum_{\textrm{irreps $w$ of $S_k$}} \chi_w(\tau)  \tilde{w} \right)$$
\end{enumerate}
\end{theorem}

\

Throughout this subsection, we will abbreviate the Young diagrams corresponding to representations of $S_n$ as indicated in Table \ref{notations}. 

\

\underline{A systematic lifting procedure for symmetric groups}

We start with a description of a systematic lifting process required for an effective implementation of Theorem \ref{ecifpp} to chains of symmetric groups. 

\

(a) Let $w_1$ and $w_2$ be two partitions of an integer $k$, represented as Young diagrams. The partial order $<$ on irreducible representations of $S_k$ can be obtained from comparing the partitions of $k$ by the number of boxes in their young diagrams below the first row. That is, set $w_1 < w_2$ if $w_1$ contains fewer boxes below the first row than $w_2$. \\
(b) We begin with an irreducible representation $w$ of $S_k$. Extend the first row of $u$ by attaching $n-k$ boxes to the first row. The result is a representation $w'$ of $S_n$ \\
(c) Observe that $\Res^{n-k}(w')$ contains exactly one copy of $w$ in its decomposition, corresponding to removing the boxes in the first row. \\
(d) Observe that all of the other representations of $S_k$ that appear $\Res^{n-k}(u')$ are $< w$ with respect to the ordering in (a) \\
(e) We therefore can apply an inductive process using $<$.

\

\begin{adjustwidth}{1cm}{}
\textbf{Remark:} It is interesting to note that the procedure we describe in (b) is the same procedure that is used in \cite{rs2} to construct the FI-module generated by a given irreducible representation. Let $w$ be a representation of $S_k$, and in the notation of \cite{rs2}, let $V(w)_n$ be the representation obtained by adding $n-k$ boxes to the top row of $w$. The combinatorial procedure of restricting $w$ to $S_k$ involves polynomial choices, and it is therefore the coefficients which describe the restriction to $S_k$ in terms of the irreducible representations of $S_k$ are polynomials in $n$. It follows that the characters $V(w)_n$ are polynomials in $n$. Exploiting this fact, the authors of \cite{rs2} show that polynomial stability is a feature of finitely generated $\FI$ modules. In our paper, the usefulness of this polynomial property is simply the observation that the coefficients that appear in our lifts are polynomials in $n$, and thus, in particular, we may simultaneously lift a representation to $S_n$ for every $n \ge k$. 
\end{adjustwidth}
\

\underline{An example of a lifting procedure:} 

\

Consider the irreducible representation $p_5 = {\tiny \yng(3,2)} \in R(S_5)$. We first attach on $n-5$ boxes to the top row of $p_5$ to get a Young diagram $p_n$, say the representation $p_9 = {\tiny \yng(7,2)}$ for $S_9$. By restricting this representation, we obtain one copy of $p_5$ as well as other representations over $S_5$: \\ ($\bullet$) There are $\choose{n-5}{1}$ ways of obtaining $v_5 = {\tiny \yng(4,1)}$ (since we may choose one step out of $n-5$ steps to remove a box from the second row) \\ ($\bullet$) And similarly, there are $\choose{n-5}{2}$ ways of of getting the trivial representation $t_5$. Thus, $p_5$ lifts to 

\begin{equation*}
\begin{split}
\tilde{p}_5 &= p_n - \choose{n-5}{1} \tilde{v}_5 - \choose{n-5}{2} \tilde{t}_5 \\
&= p_n - (n-5) (v_n-(n-5)t_n) - \frac{(n-5)(n-6)}{2} t_n \\
&= p_n - (n-5)v_n + \frac{1}{2}(n-5)(n-4)t_n
\end{split}
\end{equation*}

where we have used that $\tilde{v}_5 = v_n - (n-5)t_n$, which can be checked easily since $\Res(v_n) = t_{n-1}+v_{n-1}$.

\

Note that there are also other ways to obtain the lift of $p_5$. The lift of $v_5 \in R(S_5)$ to $S_n$ is given by
\begin{equation*}
    \tilde{v}_5 = v_n - (n-5) \tilde{t}_5
\end{equation*}
where $\tilde{t}_5$ is the unique lift of the trivial character $t_5$. The lift of $p_5 \in R(S_5)$ to $S_n$ can be obtained as:
\begin{equation*}
    \tilde{p}_5 = p_n - (n-5) \tilde{v}_5 - {\choose{n-5}{2}} \tilde{t}_5 = p_n - (n-5) v_n + \frac{1}{2}(n-5)(n-4) t_n
\end{equation*}

It should be clear from this example that we can simultaneously lift a representation of $S_k$ to every $S_n$. 

\

\underline{Formula for $\delta_{e}$}

\

Let us turn our focus to computing the character column at the identity $\delta_{e} = \C[S_n] = \sum_u \dim(u) u \in R(S_n)$. Observe that $e \in S_1$, and there is only one representation of $S_1$, the trivial representation $t_1 \in R(S_1)$, which lifts to $t \in R(S_n)$. (We will frequently commit these kinds of notational abuse, where we denote representations $w_n \in R(S_n)$ as simply $w$.) Since $\chi_t(e) = 1$, we therefore conclude from the above formula that:

$$\delta_{e} = X(X-1)...(X-(n-1)+1)t = X(X-1)...(X-(n-2))t$$

\begin{adjustwidth}{1cm}{}
\textbf{Remark:} If we think of $X = \Ind(t)$ and think of the above as tensor products, we may notice that tensoring with the trivial representation does nothing, and so one might rewrite the above as:  $\delta_{e} = X(X-1)...(X-(n-1)+1)t = X(X-1)...(X-(n-2))$. However, one should be careful to not interpret the above as a matrix. It is a vector. To avoid confusion, we will therefore stick to the first convention. 
\end{adjustwidth}

\

\underline{Formula for $\delta_{(12)}$}

\ 

$(12) \in S_2 \inj S_n$, and $S_2$ has two representations $t_2$ and $s_2$. $t_2$ lifts to $t$ and $s_2$ lifts to $s$. Moreover, $\chi_{t_2}((12)) = 1$ and $\chi_{s_2}((12)) = -1$. We therefore see that the character column of $S_n$ that corresponds to the conjugacy class of $(12)$ is: 

\ 

$$\delta_{(12)} = X(X-1)...(X-(n-3))(t-s)$$  

\ 

\underline{Formula for $\delta_{(123)}$}

\ 

$(123) \in S_3 \inj S_n$, and $S_3$ has three representations $t_3$, $s_3$, and $v_3$. $t_3$ lifts to $t$, $s_3$ lifts to $s$, and $v_3$ lifts to $v - (n-3)t$ (since $\Res(v_n) = t_{n-1}+v_{n-1}$). Observe that to compute the lift of $v_3$, we made use of our systematic lifting procedure. Observe that for representations of $S_3$, $\chi_{t_3}(123) = \chi_{s_3}(123) = 1$, and $\chi_{v_3}(123) = -1$. Thus,

\

\begin{equation*}
\begin{split}
\delta_{(123)} &= X(X-1)...(X-(n-4))(t+ s - (v-(n-3)t)) \\
&= X(X-1)...(X-(n-4))((n-2)t+s-v) \\
\end{split}
\end{equation*}

\

\underline{Formula for $\delta_{(12345)}$} 

\

Rather than going systematically, we thought it might be more instructive to provide one example of a ``slightly larger conjugacy class''. $(12345) \in S_5 \inj S_n$. The representations of $S_5$ and their respective abbreviations are given as follows:
\begin{equation*}
    \begin{array}{ccccccc}
         {\tiny \yng(5)} & {\tiny \yng(4,1)} & {\tiny \yng(3,2)} & {\tiny \yng(3,1,1)} & {\tiny \yng(2,2,1)} & {\tiny \yng(2,1,1,1)} & {\tiny \yng(1,1,1,1,1)}  \\
         t_5 & v_5 & p_5 & \wedge^2_5 & s_5p_5 & s_5v_5 & s_5
    \end{array}
\end{equation*}
For each of these representations, we get a representation of $S_n$ by stacking $n-5$ boxes to the top row of the Young diagram, and we will abuse notation a little by abbreviating representations of $S_n$, see Table \ref{notations} for further details.

\ 

Our first step is to compute lifts. Define $R(S_5)_{\ge 0} := \{ \textrm{irreps $w$ of $S_5$} \sep \chi_w(12) \ge 0\}$. For $w \in R(S_5)_{\ge 0}$, we will lift $w$ via our systematic lifting procedure. For any other representation $w'$ over $S_5$, we may observe that $s_5 w' \in R(S_5)_{\ge 0}$ and $\tilde{w'} = s_5 \cdot \tilde{(s_5 w')}$. We then obtain the lifts of representations as shown in Table \ref{table:lifts_s5}. 

\

\begin{table}[!htb]
    \centering
\begin{tabular}{ |p{3cm}||p{6 cm}| }
 \hline
 \textbf{Representation} &  \textbf{Lift} \\
  \hline
$t_5$ & $t$ \\
\hline
$v_5$ & $v - (n-5)t$ \\
\hline
$p_5$ & $p - (n-5)v + \frac{1}{2}(n-5)(n-4)t$ \\
\hline
$\wedge^2_5$ &  $\wedge^2 - (n-5)v + \frac{1}{2}(n-5)(n-4)t$ \\
\hline
$s_5p_5$ & $sp - (n-5)sv + \frac{1}{2}(n-5)(n-4)s$ \\
\hline
$s_5v_5$ & $sv - (n-5)s$ \\
\hline
$s_5$ & $s$ \\
\hline
\end{tabular}
\caption{Table of lifts for representations of $S_5$}
\label{table:lifts_s5}
\end{table}

\ 

From the character table for $S_5$, we know that $\chi_{t_5}(12345) = \chi_{s_5}(12345) = \chi_{\wedge^2_5}(12345) = 1$, $\chi_{p_5}(12345) = \chi_{s_5 p_5}(12345) = 0$, $\chi_{v_5}(12345) = \chi_{s_5 v_5}(12345) = -1$. Therefore, 

\begin{equation*}
\begin{split}
\sum_{\textrm{irreps $w$ of $S_5$}} \chi_w(12345) \tilde{w} &= t+s+\wedge^2 - (n-5)v + \frac{1}{2}(n-5)(n-4)t -(v - (n-5)t)-(sv - (n-5)s) \\ 
&= \frac{1}{2}(n-3)(n-4)t - (n-4)v + \wedge^2 - sv + (n-4)s \\ 
\end{split}
\end{equation*}

\

Therefore:

$$\delta_{(12345)} = X(X-1)...(X-(n-6))(\frac{1}{2}(n-3)(n-4)t - (n-4)v + \wedge^2 - sv + (n-4)s)$$

\

\underline{Table of formulae for character columns coming from permutations in $S_k$ for $k \le 5$}

\

Table \ref{table:formula} summarizes what have done, and provides the formula for $\delta_\tau$ for all permutations $\tau \in S_k$ for $k \le 5$.

\

\begin{table}[!htb]
\centering
\begin{tabular}{ |p{3cm}||p{11 cm}| }
 \hline
 \textbf{Permutation $\tau$} &  \textbf{Formula for $\delta_\tau$} (for every $n > k$) \\
  \hline
$e$ & $X(X-1)...(X-(n-2))t$ \newline \\ 
\hline
$(12)$ & $X(X-1)...(X-(n-3))(t-s)$ \newline \\
\hline
$(123)$ & $X(X-1)...(X-(n-4))((n-2)t+s-v)$ \newline \\
\hline
$(12)(34)$ &  $X(X-1)...(X-(n-5))((n^2 - 5n+5)t - (2n-7)v + 2p -sv + (n-4)s)$ \newline \\
\hline
$(1234)$ & $X(X-1)...(X-(n-5))((n-3)(t-s) - (v-sv))$ \newline \\
\hline
$(123)(45)$ & $X(X-1)...(X-(n-6))(\frac{(n-3)(n-4)}{2}(t-s) -(n-4)(v-sv) +(p-sp))$ \newline \\
\hline
$(12345)$ & $X(X-1)...(X-(n-6))(\frac{(n-3)(n-4)}{2}t - (n-4)v + \wedge^2 - sv + (n-4)s)$ \newline\\
\hline
\end{tabular}
\caption{Table of formulae for character columns obtained from permutations in $S_k \leq S_5$}
\label{table:formula}
\end{table}

\

\underline{Computing the character column $\delta_{(123)}$ of $S_6$}

In this example, we explicitly compute the character column $\delta_{(123)}$ of $S_6$. The representations of $S_6$ and the corresponding abbreviations are given as follows.
\begin{equation*}
    \begin{array}{ccccccccccc}
       {\tiny \yng(6)} & {\tiny \yng(5,1)} & {\tiny \yng(4,2)} & {\tiny \yng(4,1,1)} & {\tiny \yng(3,3)} & {\tiny \yng(3,2,1)} & {\tiny \yng(2,2,2)} & {\tiny \yng(3,1,1,1)} & {\tiny \yng(2,2,1,1)} & {\tiny \yng(2,1,1,1,1)} & {\tiny \yng(1,1,1,1,1,1)} \\
    t_6 & v_6 & p_6 & \wedge^2_6 & b_6 & r_6 & s_6 b_6 & s_6 \wedge^2_6 & s_6 p_6 & s_6 v_6 & s_6
    \end{array}
\end{equation*}

\

Consider the McKay graph $\M(S_6, X)$, which describes the irreducible decomposition of tensor products of representations of $S_6$. The vertices of the graph $\M(G, V)$ are the irreducible representations of $G$, and the number of edges from $U$ to $W$ are $\inner{U \x V}{W}$. The McKay graph serves as a visual aid for constructing the matrix $X = \Ind \Res$. In addition, counting paths and estimating numbers of paths on the McKay graph might play a crucial role in combinatorial or probabilistic extensions of this paper.
\

We may easily construct the McKay graph $\M(S_6, X)$ using the Branching Rule. For example, 

$$\Ind \Res({\tiny \yng(4,2))} = \Ind ({\tiny \yng(3,2)} + {\tiny \yng(4,1)}) = ({\tiny \yng(4,2)} + {\tiny \yng(3,3)} + {\tiny \yng(3,2,1)})+ ({\tiny \yng(5,1)} + {\tiny \yng(4,2)}+{\tiny \yng(4,1,1)})$$ and therefore, $X p_6 = 2p_6 + v_6 + b_6 + \wedge^2_6 + r_6$. The reader may verify that we obtain the following graph as shown in Figure \ref{figure:mckay}, where the subscripts indicating that these are representations over $S_6$ are abbreviated.

\

\begin{figure}[!htb]
    \centering
    \begin{tikzpicture}[
            > = stealth, % arrow head style
            shorten > = 1pt, % don't touch arrow head to node
            auto,
            node distance = 2.75cm, % distance between nodes
            semithick % line style
        ]

        \tikzstyle{every state}=[
            draw = black,
            thick,
            fill = white,
            minimum size = 80mm
        ]

        \node (t_6) [circle, draw = black]{$t$};
        \node (v_6) [circle, draw=black, above right of=t_6] {$v$};
        \node (p_6) [circle, draw=black, below right of=t_6] {$p$};
        \node (wv_6) [circle, draw=black, right of=v_6] {$\wedge^2$};
        \node (b_6) [circle, draw=black, right of=p_6] {$b$};
        \node (r_6) [circle, draw=black, above right of=b_6] {$r$};
        \node (sb_6)[circle, draw=black, below right of=r_6] {$s b$};
        \node (swv_6) [circle, draw=black, above right of=r_6] {$s \wedge^2$};
        \node (sp_6) [circle, draw=black, right of=sb_6] {$s p$};
        \node (sv_6) [circle, draw=black, right of=swv_6] {$s v$};
        \node (s_6) [circle, draw=black, below right of=sv_6] {$s$};

        \path[-] (t_6) edge [loop above] node {1} (t_6);
        \path[-] (t_6) edge node {1} (v_6);
        \path[-] (v_6) edge [loop above] node {2} (v_6);
        \path[-] (v_6) edge node {1} (p_6);
        \path[-] (v_6) edge node {1} (wv_6);
        \path[-] (p_6) edge [loop below] node {2} (p_6);
        \path[-] (p_6) edge node {1} (r_6);
        \path[-] (p_6) edge node {1} (b_6);
        \path[-] (p_6) edge node {1} (wv_6);
        \path[-] (wv_6) edge [loop above] node {2} (wv_6);
        \path[-] (wv_6) edge node {1} (r_6);
        \path[-] (wv_6) edge node {1} (swv_6);
        \path[-] (b_6) edge node {1} (r_6);
        \path[-] (b_6) edge [loop below] node {1} (b_6);
        \path[-] (r_6) edge [loop above] node {3} (r_6);
        \path[-] (r_6) edge node {1} (sb_6);
        \path[-] (r_6) edge node {1} (sp_6);
        \path[-] (r_6) edge node {1} (swv_6);
        \path[-] (sb_6) edge [loop below] node {1} (sb_6);
        \path[-] (sb_6) edge node {1} (sp_6);
        \path[-] (swv_6) edge [loop above] node {2} (swv_6);
        \path[-] (swv_6) edge node {1} (sv_6);
        \path[-] (swv_6) edge node {1} (sp_6);
        \path[-] (sp_6) edge [loop below] node {2} (sp_6);
        \path[-] (sp_6) edge node {1} (sv_6);
        \path[-] (sv_6) edge [loop above] node {2} (sv_6);
        \path[-] (sv_6) edge node {1} (s_6);
        \path[-] (s_6) edge [loop above] node {1} (s_6);
        
        %\path[-] (rho0) edge node {} (rho1);
    \end{tikzpicture}
    \caption{McKay graph $\M(S_6, X)$}
    \label{figure:mckay}
\end{figure}

$X$ is then the following transition matrix of the Mckay graph: 
\begin{equation*}
{\small
    X = 
    \begin{pmatrix}
    1 & 1 & 0 & 0 & 0 & 0 & 0 & 0 & 0 & 0 & 0 \\
    1 & 2 & 1 & 1 & 0 & 0 & 0 & 0 & 0 & 0 & 0 \\
    0 & 1 & 2 & 1 & 1 & 1 & 0 & 0 & 0 & 0 & 0 \\
    0 & 1 & 1 & 2 & 0 & 1 & 0 & 1 & 0 & 0 & 0 \\
    0 & 0 & 1 & 0 & 1 & 1 & 0 & 0 & 0 & 0 & 0 \\
    0 & 0 & 1 & 1 & 1 & 3 & 1 & 1 & 1 & 0 & 0 \\
    0 & 0 & 0 & 0 & 0 & 1 & 1 & 0 & 1 & 0 & 0 \\
    0 & 0 & 0 & 1 & 0 & 1 & 0 & 2 & 1 & 1 & 0 \\
    0 & 0 & 0 & 0 & 0 & 1 & 1 & 1 & 2 & 1 & 0 \\
    0 & 0 & 0 & 0 & 0 & 0 & 0 & 1 & 1 & 2 & 1 \\
    0 & 0 & 0 & 0 & 0 & 0 & 0 & 0 & 0 & 1 & 1 \\
    \end{pmatrix}}
\end{equation*}

\

Theorem \ref{ecifpp} thus implies that:

\

\begin{equation*}
\begin{split}
\delta_{(123)} &= X(X-1)....(X-(n-4))((n-2)t+s-v) \\
&= X(X-1)....(X-(n-4)) \times \begin{pmatrix} n-2 & -1 & 0 & 0 & 0 & 0 & 0 & 0 & 0 & 0 & 1 \end{pmatrix}^T \\ 
\end{split}
\end{equation*}

\

Substituting the matrix $X$, we derive the character column of the conjugacy class $(123)$ for $S_6$:

\ 

$$\delta_{(123)} = \begin{pmatrix} 1 & 2 & 0 & 1 & -1 & -2 & -1 & 1 & 0 & 2 & 1 \end{pmatrix}^T$$

\

\subsection{Chains of wreath products}

Let $H$ be a finite group. We give a quick review of the representation theory of wreath products $H^n \rtimes S_n$. We will introduce notation which, in our opinion, is both conceptually clear and computationally ideal. We refer the reader to \cite{cw2}, \cite{cw3}, \cite{cw1}, and \cite{lr1} for more details. A representation of the wreath product $H^n \sdp S_n$ is denoted by a $2 \times d$ array
\begin{equation*}
    \begin{pmatrix} U_1 & U_2 & ... & U_d \\ \lambda_1 & \lambda_2 & ... & \lambda_d \end{pmatrix} := \Ind_{H^n \sdp (S_{k_1} \times \cdots \times S_{k_d})}^{H^n \sdp S_n} (U_1^{\otimes k_1} \otimes \lambda_1) \otimes \cdots \otimes (U_d^{\otimes k_d} \otimes \lambda_d),
\end{equation*}
where $U_i$ is an irreducible representation of $H$, and $\lambda_i$ is a representation of $S_{k_i}$. The specific order of the irreducible representations $\{U_i\}_{i=1}^n$ of $H$ is fixed. The irreducible representations of the wreath product $H^n \sdp S_n$ can be denoted by:

$$\left\{ \begin{pmatrix} U_1 & U_2 & ... & U_d \\ \lambda_1 & \lambda_2 & ... & \lambda_d \end{pmatrix} \sep \{U_i\}_{i=1}^d \textrm{ are distinct irreps of $H$}, \lambda_i \textrm{ is an irrep of $S_{k_i}$ for each $i$} \right\}$$.

\

The action of $H^n \sdp S_n$ on the array $\begin{pmatrix} U_1 & U_2 & ... & U_d \\ \lambda_1 & \lambda_2 & ... & \lambda_d \end{pmatrix}$ is as follows: $H^n$ and $S_{k_1} \times ... \times S_{k_d}$ both act on $U_1 \x U_2 \x .. \x U_d$; $H^n$ acts on each factor and $S_{k_1} \times ... \times S_{k_d}$ permutes the factors according to the representation $\lambda_1 \x ... \x \lambda_d$ (Note that $\lambda_i$ is a representation of $S_{k_i}$ for each $i$).

\

The branching rule, which plays a crucial role for us, is easy to state using our array notation:

\

$$\Ind \left(\begin{pmatrix} U_1 & U_2 & ... & U_d \\ \lambda_1 & \lambda_2 & ... & \lambda_d \end{pmatrix} \right) = \sum_{i=1}^d \dim(U_i) \cdot \begin{pmatrix} U_1 & ... & U_i & ... & U_d \\ \lambda_1 & ... & \Ind(\lambda_i) & ... & \lambda_d \end{pmatrix}$$

$$\Res \left( \begin{pmatrix} U_1 & U_2 & ... & U_d \\ \lambda_1 & \lambda_2 & ... & \lambda_d \end{pmatrix} \right) = \sum_{i=1}^d \dim(U_i) \cdot \begin{pmatrix} U_1 & ... & U_i & ... & U_d \\ \lambda_1 & ... & \Res(\lambda_i) & ... & \lambda_d \end{pmatrix}$$

\

We recall that the irreducible representations of $S_n$ are represented by Young diagrams. The branching rules for induction and restriction are given by adding or removing a box from the given diagram.

\

The reader should take note that there is a subtle point in the above rule: namely that we may restrict a Young diagram with only $1$ box to a Young diagram with $0$ boxes, and thereby removing the corresponding column. Likewise, when inducing a representation whose Young diagram does not have a full set of columns, we may add additional columns of boxes where necessary. Adopting these notations, one can prove the following two facts about chains of wreath products, the proof of which we omit in this paper.
\begin{proposition} \label{cwpsc}
Let $H$ be any finite group. 
\begin{enumerate}
    \item The chain $\{H^n \sdp S_n\}_{n \in \N}$ is a surjective chain. 
    \item The chain $\{H^n \sdp S_n\}_{n \in \N}$ satisfies the Heisenberg algebra property with scaling $|H|$.
\end{enumerate}
\end{proposition}

For any group $H$, it is well known that the chain of wreath products $\{H^n \sdp S_n\}_{n \in \N}$ satisfies the Heisenberg algebra property with scaling $|H|$. It should be noted that the only known examples of chains of groups that satisfy the Heisenberg algebra property are chains of wreath products. It is proved in \cite{dg4} that if $M = 1$ or $M$ is prime, then the only such chains are chains of wreath products. It is further conjectured in \cite{dg4} that the only such chains are chains of wreath products for $M$ arbitrary. 

The statement of the main results can be applied to chains of wreath products as follows:
\begin{theorem}
Let $H$ be any finite group. Let $X = \Ind \Res = \Ind(t) \otimes$ be the $\Ind \Res$ operator on the chain of wreath products $\{H^n \rtimes S_n\}$.
\begin{enumerate}
    \item The polynomial $f_l$ is given by
    $$f_l(X) = X(X-|H|)(X-2|H|)...(X-(l-1)|H|)$$
    \item The non-identity character values of $\Ind t$ in the wreath product $H^n \rtimes S_n$ is the set $\{0, |H|, \cdots, (n-1)|H|\}$.
    \item Let $\alpha \in H^k \sdp S_k \inj H^{k+l} \sdp S_{k+l}$, and let $\delta_\alpha \in R(H^{k+l} \sdp S_{k+l})$ be the character column of $\alpha$, that is: $\delta_\alpha = \sum_{\textrm{irreps $u$ of $H^{k+l} \sdp S_{k+l}$}} \chi_u(\alpha) u$. For each representation $w$ of $G_k$, let $\tilde{w}$ denote a lift to $G_{k+l}$, i.e an element of $R(G_{k+l})$ such that $\Res^l(\tilde{w}) = w$.Then: 
$$\delta_\alpha = \Ind \Res ( \Ind \Res - |H|)(\Ind \Res-2|H|)...(\Ind \Res-(l-1)|H|) \left(\sum_{\textrm{irreps $w$ of $H^k \sdp S_k$}} \chi_w(\alpha) \cdot \tilde{w} \right)$$
\end{enumerate}
\end{theorem}

\underline{A systematic lifting procedure for wreath product representations} \label{slp2}

\

As shown for chains of symmetric groups, we begin with a blueprint of the systematic lifting of representations for chains of wreath products.

\

(a) For a Young diagram $\lambda$ let $B(\lambda)$ denote the number of boxes of $\lambda$ and let $B'(\lambda)$ denote the number of boxes of $\lambda$ below the top row. \\
(b) Define the partial ordering $<$ on the irreducible representations of $H^k \sdp S_k$ by the following rule: \\

\begin{equation*}
\begin{split}
&\begin{pmatrix} U_1 & U_2 & ... & U_d \\ \lambda_1 & \lambda_2 & ... & \lambda_d \end{pmatrix} < \begin{pmatrix} U_1 & U_2 & ... & U_d \\ \mu_1 & \mu_2 & ... & \mu_d \end{pmatrix} \textrm{ if either of the following holds:} \\
& \textrm{(1) $B'(\lambda_1) < B'(\mu_1)$ and $B(\lambda_i) \le B(\mu_i)$ for all $i \ge 2$.} \\ 
& \textrm{(2) $B'(\lambda_1) \le B'(\mu_1)$, $B(\lambda_i) \le B(\mu_i)$ for all $i \ge 2$, and $B(\lambda_i) < B(\mu_i)$ for at least one $i \ge 2$.} \\
\end{split}
\end{equation*}

where a subtle point is that we allow $\mu_i = 0$ or $\lambda_i = 0$ (the Young diagram with $0$ boxes), since the branching rule of representations of symmetric groups allows us to restrict a Young diagram with only one box to obtain the Young diagram with $0$ boxes, which corresponds to removing the corresponding column. \\

(c) We begin with a representation $w = \begin{pmatrix} U_1 & U_2 & ... & U_d \\ \mu_1 & \mu_2 & ... & \mu_d \end{pmatrix}$ of $H^k \sdp S_k$. Extend the first row of $\mu_1$ by attaching $n-k$ boxes to the first row, and call the resulting Young diagram $\mu_1'$. Consider the element $$w' = \frac{1}{\dim(U_1)^{n-k}} \begin{pmatrix} U_1 & U_2 & ... & U_d \\ \mu_1' & \mu_2 & ... & \mu_d \end{pmatrix} \in R(H^n \sdp S_n)$$ \\
(d) Observe that $\Res^{n-k}(w')$ contains exactly one copy of $w$ in its decomposition, corresponding to removing the boxes in the first row of $\mu_1$. \\
(e) Observe that all of the other representations of $S_k$ that appear $\Res^{n-k}(u')$ are $< w$ with respect to the ordering in (a) \\
(f) We therefore can apply an inductive process using $<$.

\

\underline{An example to illlustrate the procedure above:}

\

Let $1$ and $-1$ be the representations of $\quotient{\Z}{2\Z}$. Denote by $t_i$ and $s_i$ the trivial and the sign representations of $S_i$. We will lift the representation $\begin{pmatrix} 1 & -1 \\ t_1 & t_1 \end{pmatrix} \in R(\left(\quotient{\Z}{2\Z}\right)^2 \sdp S_2)$ to $R(\left(\quotient{\Z}{2\Z}\right)^n \sdp S_n)$ via our systematic lifting procedure. The first step is to consider the representation $\begin{pmatrix} 1 & -1 \\ t_{n-1} & t_1 \end{pmatrix} \in R(\left(\quotient{\Z}{2\Z}\right)^n \sdp S_n)$, and observe that:

\

$$\Res^{n-2} \begin{pmatrix} 1 & -1 \\ t_{n-1} & t \end{pmatrix} = \begin{pmatrix} 1 & -1 \\ t_1 & t_1 \end{pmatrix} + (n-2) \begin{pmatrix} 1 \\ t_2 \end{pmatrix}$$

since there are $\choose{n-2}{1} = (n-2)$ ways to remove a box from the second column (one way for each step). 

\

Since $\begin{pmatrix} 1 \\ t_2 \end{pmatrix}$ is the trivial representation, which lifts to the trivial representation $\begin{pmatrix} 1 \\ t_n \end{pmatrix}$, we conclude that: 

\

$$\begin{pmatrix} 1 & -1 \\ t_1 & t_1 \end{pmatrix} \textrm{ lifts to } \begin{pmatrix} 1 & -1 \\ t_{n-1} & t_1 \end{pmatrix} - (n-2) \begin{pmatrix} 1 \\ t_n \end{pmatrix}$$ 

\

\underline{Character columns of $(\Z/2\Z)^n \sdp S_n$}

\

We will let $\pm 1$ denote both the elements and the representations of the group $ \Z/2\Z$, but it will be clear when we are referring to which. As before, denote by $t_i$ and $s_i$ the trivial and the sign representations of $S_i$. $\Z/2\Z \sdp S_1 \isom \Z/2\Z$ has two representations, $\begin{pmatrix} 1 \\ t_1 \end{pmatrix}$ and $\begin{pmatrix} -1 \\ t_1 \end{pmatrix}$. The representation $\begin{pmatrix} 1 \\ t_1 \end{pmatrix}$ is the trivial representation which lifts to the trivial representation  $\begin{pmatrix} 1 \\ t_n \end{pmatrix} \in R((\Z/2\Z)^n \sdp S_n)$. $\begin{pmatrix} -1 \\ t_1 \end{pmatrix}$ lifts to $\begin{pmatrix} -1 \\ t_n \end{pmatrix}$. 

\

The characters of representations of $\Z/2\Z \sdp S_1$ are: $$\chi_{\begin{pmatrix} 1 \\ t_1 \end{pmatrix}}((1), e) = \chi_{\begin{pmatrix} 1 \\ t_1 \end{pmatrix}}((-1), e) = \chi_{\begin{pmatrix} -1 \\ t_1 \end{pmatrix}}((1), e) = 1$$ $$\chi_{\begin{pmatrix} -1 \\ t_1 \end{pmatrix}}((-1), e) = -1$$

\

Therefore,

\

$$\delta_e = X(X-2)(X-4)...(X- 2(n-2)) \left( \begin{pmatrix} 1 \\ t_n \end{pmatrix} + \begin{pmatrix} -1 \\ t_n \end{pmatrix} \right)$$
$$\delta_{((-1,1,1,...,1),e )} = X(X-2)(X-4)...(X- 2(n-2)) \left(\begin{pmatrix} 1 \\ t_n \end{pmatrix} - \begin{pmatrix} -1 \\ t_n \end{pmatrix} \right)$$

\

To obtain character columns $\delta_{((-1,-1,1,...,1),e )}$, $\delta_{((1,1,1,...,1),(12) )}$, and $\delta_{((-1,1,1,...,1),(12) )}$, we refer to the character table of $(\Z/2\Z)^2 \sdp S_2$, given as in Table \ref{table:wreath}.

\

\begin{table}[!htb]
\centering
\begin{tabular}{ |p{2 cm}|p{2 cm}|p{2 cm}|p{2 cm}|p{2 cm}|p{2 cm}|}
 \hline
  &  $e$ & $((-1,1), e)$ & $((-1,-1), e)$ & $((1,1), (12))$ & $((-1,1), (12))$ \\
  \hline 
  \break $\begin{pmatrix} 1 \\ t_2 \end{pmatrix}$ \newline &  \break 1 \newline & \break 1 \newline & \break 1 \newline &  \break 1 \newline & \break 1 \newline \\
  \hline
   \break $\begin{pmatrix} 1 & -1 \\ t_1 & t_1 \end{pmatrix}$ \newline &  \break 2 \newline & \break 0 \newline & \break -2 \newline &  \break 0 \newline & \break 0 \newline \\
  \hline
  \break $\begin{pmatrix} 1 \\ s_2 \end{pmatrix}$ \newline &  \break 1 \newline & \break 1 \newline & \break 1 \newline &  \break -1 \newline & \break -1 \newline \\
  \hline
  \break $\begin{pmatrix} -1 \\ t_2 \end{pmatrix}$ \newline &  \break 1 \newline & \break -1 \newline & \break 1 \newline &  \break 1 \newline & \break -1 \newline \\
  \hline
  \break $\begin{pmatrix} -1 \\ s_2 \end{pmatrix}$ \newline &  \break 1 \newline & \break -1 \newline & \break 1 \newline &  \break -1 \newline & \break 1 \newline \\
  \hline 
\end{tabular}
\caption{Character table for $(\Z/2\Z)^2 \sdp S_2$}
\label{table:wreath}
\end{table}

\

Given any one dimensional representation $U$ of $H$, the branching rule of representations of wreath products imply that $\begin{pmatrix} U \\ t_k \end{pmatrix} \in R(H^k \sdp S_k)$ lifts to $\begin{pmatrix} U \\ t_n \end{pmatrix} \in R(H^n \sdp S_n)$ and $\begin{pmatrix} U \\ s_k \end{pmatrix} \in R(H^k \sdp S_k)$ lifts to $\begin{pmatrix} U \\ s_n \end{pmatrix} \in R(H^n \sdp S_n)$. Therefore, all lifts except for the lift of $\begin{pmatrix} 1 & -1 \\ t_1 & t_1 \end{pmatrix}$ can be obtained. This can be done via the proposed systematic lifting procedure, and we obtain that $\begin{pmatrix} 1 & -1 \\ t_1 & t_1 \end{pmatrix}$ lifts to $\begin{pmatrix} 1 & -1 \\ t_{n-1} & t_1 \end{pmatrix} - (n-2) \begin{pmatrix} 1 \\ t_n \end{pmatrix}$. 

\

Therefore, we have: 

\

\begin{equation*}
\begin{split}
\delta_{((-1,-1,1,...,1),e )} = X(X-2)...(X-2(n-3)) \cdot \bigg[ & (2n-3) \begin{pmatrix} 1 \\ t_n \end{pmatrix} -2 \begin{pmatrix} 1 & -1 \\ t_{n-1} & t_1 \end{pmatrix} + \begin{pmatrix} 1 \\ s_n \end{pmatrix} + \begin{pmatrix} -1 \\ t_n \end{pmatrix}  + \begin{pmatrix} -1 \\ s_n \end{pmatrix} \bigg] \\
\end{split}
\end{equation*}

\begin{equation*}
\begin{split}
\delta_{((1,1,1,...,1),(12) )} = X(X-2)...(X-2(n-3)) \cdot \bigg[ & \begin{pmatrix} 1 \\ t_n \end{pmatrix} - \begin{pmatrix} 1 \\ s_n \end{pmatrix}  + \begin{pmatrix} -1 \\ t_n \end{pmatrix} - \begin{pmatrix} -1 \\ s_n \end{pmatrix} \bigg] \\
\end{split}
\end{equation*}

\begin{equation*}
\begin{split}
\delta_{((-1,1,1,...,1),(12) )} = X(X-2)...(X-2(n-3)) \cdot \bigg[ & \begin{pmatrix} 1 \\ t_n \end{pmatrix} - \begin{pmatrix} 1 \\ s_n \end{pmatrix}  - \begin{pmatrix} -1 \\ t_n \end{pmatrix} + \begin{pmatrix} -1 \\ s_n \end{pmatrix} \bigg] \\
\end{split}
\end{equation*}

\

\underline{Character columns of $(\Z/2\Z)^3 \rtimes S_3$}

\

We compute character columns of $(\Z/2\Z)^3 \rtimes S_3$ for $\delta_e, \delta_{((-1, 1, 1), e)}$, $\delta_{((-1,-1,1),e)}$, $\delta_{((1,1,1),e)}$, and $\delta_{((-1,1,1),(12))}$ as an example. Using the branching rule, we compute the matrix $X$ can be written as follows with respect to the following basis of irreducible representations:

$$\left\{ \begin{pmatrix} 1 \\ t_3 \end{pmatrix}, \begin{pmatrix} 1 \\ s_3 \end{pmatrix}, \begin{pmatrix} 1 \\ v_3 \end{pmatrix}, \begin{pmatrix} 1 & -1 \\ t_2 & t_1 \end{pmatrix}, \begin{pmatrix} 1 & -1 \\ s_2 & t_1 \end{pmatrix}, \begin{pmatrix} 1 & -1 \\ t_1 & t_2 \end{pmatrix}, \begin{pmatrix} 1 & -1 \\ t_1 & s_2 \end{pmatrix}, \begin{pmatrix} -1 \\ t_3 \end{pmatrix}, \begin{pmatrix} -1 \\ s_3 \end{pmatrix}, \begin{pmatrix} -1 \\ v_3 \end{pmatrix} \right\}$$
\begin{equation*}
X = 
    \begin{pmatrix}
    1 & 0 & 1 & 1 & 0 & 0 & 0 & 0 & 0 & 0 \\
    0 & 1 & 1 & 0 & 1 & 0 & 0 & 0 & 0 & 0 \\
    1 & 1 & 2 & 1 & 1 & 0 & 0 & 0 & 0 & 0 \\
    1 & 0 & 1 & 2 & 1 & 1 & 1 & 0 & 0 & 0 \\
    0 & 1 & 1 & 1 & 2 & 1 & 1 & 0 & 0 & 0 \\
    0 & 0 & 0 & 1 & 1 & 2 & 1 & 1 & 0 & 1 \\
    0 & 0 & 0 & 1 & 1 & 1 & 2 & 0 & 1 & 1 \\
    0 & 0 & 0 & 0 & 0 & 1 & 0 & 1 & 0 & 1 \\
    0 & 0 & 0 & 0 & 0 & 0 & 1 & 0 & 1 & 1 \\
    0 & 0 & 0 & 0 & 0 & 1 & 1 & 1 & 1 & 2
    \end{pmatrix}
\end{equation*}
We then have:
\begin{align*}
    \delta_e &= X(X-2) \times \begin{pmatrix} 1 & 0 & 0 & 0 & 0 & 0 & 0 & 1 & 0 & 0 \end{pmatrix}^T \\
    &= \begin{pmatrix} 1 & 1 & 2 & 3 & 3 & 3 & 3 & 1 & 1 & 2 \end{pmatrix}^T
\end{align*}
\begin{align*}
    \delta_{((-1,1,1),e)} &= X(X-2) \begin{pmatrix} 1 & 0 & 0 & 0 & 0 & 0 & 0 & -1 & 0 & 0 \end{pmatrix}^T \\
    &= \begin{pmatrix} 1 & 1 & 2 & 1 & 1 & -1 & -1 & -1 & -1 & -2 \end{pmatrix}^T
\end{align*}
\begin{align*}
    \delta_{((-1,-1,1),e)} &= X \times \begin{pmatrix} 3 & 1 & 0 & -2 & 0 & 0 & 0 & 1 & 1 & 0 \end{pmatrix}^T \\
    &= \begin{pmatrix} 1 & 1 & 2 & -1 & -1 & -1 & -1 & 1 & 1 & 2 \end{pmatrix}^T
\end{align*}
\begin{align*}
    \delta_{((1, 1, 1),(1 2))} &= X \begin{pmatrix} 1 & -1 & 0 & 0 & 0 & 0 & 0 & 1 & -1 & 0 \end{pmatrix}^T \\
    &= \begin{pmatrix} 1 & -1 & 0 & 1 & -1 & 1 & -1 & 1 & -1 & 0 \end{pmatrix}^T
\end{align*}
\begin{align*}
    \delta_{((-1, 1, 1),(1 2))} &= X \begin{pmatrix} 1 & -1 & 0 & 0 & 0 & 0 & 0 & -1 & 1 & 0 \end{pmatrix}^T \\
    &= \begin{pmatrix} 1 & -1 & 0 & 1 & -1 & -1 & 1 & -1 & 1 & 0 \end{pmatrix}^T
\end{align*}

\
\

\section{Future directions}

\

We discuss three questions which naturally arise as a consequence of this paper:  

\

(1) As a consequence of Theorem \ref{pdtp}, we know that if a surjective chain satisfies \hyperref[*1]{Property $(*)$}, then the polynomials $\{f_l\}_{l \in \N}$ are determined by two parameters $B \in \N$, $C \in \Z$. Chains of wreath products show that it is possible to construct chains corresponding to $(B,C) = (1,c)$, for $c \in \Z$ arbitrary. A natural question to ask is: What tuples $(B,C)$ are realized by some chain of groups? Prerequisite questions are: Can $B$ vary, or is $B = 1$ forced? Also, is $C$ forced to be a positive integer, as in the case where $G_0$ is trivial? Note that in the special case of a dual tower of groups, it is conjectured in \cite{dg4} that every such chain must be a chain of wreath products, and thus $B=1$ is conjecturally forced. It is proved in \cite{dg4} that this conjecture is true whenever the scaling factor is either $1$ or is a prime. However, the property of being a dual tower of groups is a very restrictive property, and we are interested in \hyperref[*1]{Property $(*)$}, which is a far less restrictive property.

\

(2) Can we make asymptotic statements on character columns $\delta_\alpha$ by applying probabilistic techniques in conjunction with our formulae? Since McKay graphs of symmetric groups stabilize as $n \to \infty$, the matrix entries of $X = \Ind \Res$ stabilize as $n \to \infty$. Furthermore, our method fixes a conjugacy class  $\alpha$ and simultaneously produces a formula for the corresponding column in the character table of $S_n$ for every $n$. Hence, we can ask: what is the probability distribution of $\delta_\alpha$ as $n \to \infty$? It should be known that local asymptotics (i.e asymptotics for a fixed limit shape of a representation) is well understood \cite{as1}. However, there is much which is yet to be understood about the global asymptotics (i.e the asymptotic behavior of all possible limit shapes and their relations). Since our formulae produce entire character columns at once, we produce global data, and hence it is meaningful to ask whether our formulae are useful in the context of global asymptotics.

\

(3) How can we interpret our formulae as generalizations of projection operators? One can see that $X(X-1)..(X-(n-2))$ is in fact the projection operator onto $\delta_{(e)} \in R(S_n)$ and $(t-s)X(X-1)..(X-(n-2))$ is in fact the projection operator onto $\delta_{(12)} \in R(S_n)$. In general however, our operators are not projection operators onto lines, but do project certain subspaces onto lines. We believe it would be insightful to investigate this line of thought further.  

\section{Acknowledgements} 

\

We would like to thank Mikhail Khovanov and Jordan Ellenberg for their insightful comments. We would like to thank the reviewers for carefully examining the manuscript and giving constructive feedback and comments.

\bibliography{D26PROBCCR}{}
\bibliographystyle{plain}

\end{document}